\documentclass{amsart}

\usepackage{graphicx}
\usepackage{caption}
\usepackage{subcaption}
\RequirePackage{amsfonts,amssymb,amsbsy,amsmath,amscd}

\usepackage{todonotes}
\newcommand{\C}{\mathbb{C}}
\newcommand{\PP}{\mathbb{P}}
\newcommand{\R}{\mathbb{R}}
\newcommand{\Z}{\mathbb{Z}}
\newcommand{\Ff}{\mathcal{F}}
\newcommand{\Oo}{\mathcal{O}}
\newcommand{\Mm}{\mathcal{M}}
\newcommand{\Tt}{\mathcal{T}}
\newcommand{\Ss}{\mathcal{S}}
\newcommand{\Aa}{\mathcal{A}}
\newcommand{\Cc}{\mathcal{C}}
\newcommand{\eps}{\varepsilon}
\newcommand{\Res}{\mathrm{Res}}
\newcommand{\de}{\partial}
\newcommand{\proc}[1]{\medbreak\noindent{\it #1}\hspace{1ex}\ignorespaces}
\newcommand{\ep}{\noindent{\hfill $\Box$}}
\newcommand{\affilnum}[1]{\ifcase#1\or $\dagger$\or $\ddagger$\or
   $\S$\or $\P$\or $\|$\or **\or $\dagger\dagger$
   \or $\ddagger\ddagger$\else\fi}
\newcommand{\affil}[1]{\ifcase#1\or $\dagger$\or $\ddagger$\or
   $\S$\or $\P$\or $\|$\or **\or $\dagger\dagger$
   \or $\ddagger\ddagger$\else\fi}
\renewcommand{\Re}{\mathrm{Re}\,}
\renewcommand{\Im}{\mathrm{Im}\,}
\newtheorem{theo}{Theorem}[section]
\newtheorem{prop}{Proposition}[section]
\newtheorem{lemma}{Lemma}[section]
\newtheorem{cor}{Corollary}[section]
\numberwithin{equation}{section}

\begin{document}


\title[Meromorphic connections on Riemann surfaces]{A Poincar\'e-Bendixson theorem for meromorphic connections on Riemann surfaces}


\author{Marco Abate}
\address{Dipartimento di Matematica, Universit\`a di Pisa, Largo Pontecorvo 5,
56127 Pisa, Italy}
\email{abate@dm.unipi.it}

\author{Fabrizio Bianchi}
\address{Institut de Math\'ematiques de Toulouse, Universit\'e Paul Sabatier, 118 route de Narbonne, F-31062 Toulouse, France}
\email{fbianchi@math.univ-toulouse.fr}

\thanks{Both authors were partially supported by the FIRB 2012 project \emph{Differential geometry and geometric function theory.}}

\begin{abstract}
We shall prove a Poincar\'e-Bendixson theorem describing the asymptotic behavior of geodesics
for a meromorphic connection on a compact
Riemann surface. We shall also briefly discuss the case of non-compact Riemann surfaces, and
study in detail the geodesics for a holomorphic connection on a complex torus.  
\end{abstract}

\maketitle

\section{Introduction}

The main goal of this paper is to prove a Poincar\'e-Bendixson theorem describing the asymptotic behavior
of geodesics for a meromorphic connection on a compact
Riemann surface.

Roughly speaking (see Section~\ref{section_preliminary} for more details), a \emph{meromorphic connection}
on a Riemann surface $S$ is a $\C$-linear operator $\nabla\colon \mathcal{T}(S)\to\mathcal{M}^1(S)\otimes\mathcal{T}(S)$,
where $\mathcal{T}(S)$ is the space of holomorphic vector fields on~$S$ and $\mathcal{M}^1(S)$ is the space of meromorphic
1-forms on~$S$, satisfying a Leibniz rule $\nabla(\phi s)=\phi\nabla s+d\phi\otimes s$ for every $s\in\mathcal{T}(S)$ and
every meromorphic function $\phi\in\mathcal{M}(S)$.
Meromorphic connections are a classical subject, related for instance to linear differential systems (see, e.g., \cite{ilyako});
but in this paper we shall study them following a less classical point of view, introduced in~\cite{ab}.

A meromorphic connection~$\nabla$ on a Riemann surface~$S$ (again, see Section~\ref{section_preliminary} for details) can
be represented in local coordinates by a meromorphic 1-form; it turns out that the poles (and the associated residues)
of this form do not depend on the coordinates but only on the connection. 
If $p\in S$ is not a pole, as for ordinary connections in differential geometry, it is possible to use $\nabla$ for
differentiating along a direction $v\in T_pS$ any vector field $s$ defined along a curve in~$S$ tangent to~$v$, obtaining
a tangent vector $\nabla_v s\in T_pS$. In particular, if $\sigma\colon (a,b)\to S$ is a smooth curve with image contained
in the complement of the poles, it makes sense to consider the vector field $\nabla_{\sigma'}\sigma'$ along~$\sigma$; and we
shall say that $\sigma$ is a \emph{geodesic} for~$\nabla$ if $\nabla_{\sigma'}\sigma'\equiv O$.

As far as we know, geodesics for meromorphic connections in this sense were first introduced in~\cite{ab}. As shown
there, locally they behave as Riemannian geodesics of a flat metric; but their global behavior can be very different from the
global behavior of Riemannian geodesics. Thus they are an interesting object of study \textit{per se}; but they
also have dynamical applications, in particular in the theory of local discrete holomorphic dynamical systems of several
variables --- explaining why we are particularly interested in their asymptotic behavior.

One of the main open problems in local dynamics of several complex variables is the
understanding of the dynamics, in a full neighbourhood of the origin, of holomorphic germs tangent to the identity, that is
of germs of holomorphic endomorphisms of~$\C^n$ fixing the origin and with differential there
equal to the identity. In dimension one, the Leau-Fatou flower theorem (see, e.g., \cite{abasurvey} or \cite{milnor})
provides exactly such an understanding; and building on this theorem Camacho (\cite{camacho78}; see also \cite{shcherbakov1982topological}) in
1978 has been able to prove that every holomorphic germ tangent to the identity in dimension one is locally topologically conjugated
to the time-1 map of a homogeneous vector field. In other words, time-$1$ maps of homogeneous vector fields provide a complete list
of models for the local topological dynamics of one-dimensional holomorphic germs tangent to the identity.

In recent years, many authors have begun to study the local dynamics of germs tangent to the identity in several complex variables;
see, e.g., \'Ecalle \cite{ecalle1, ecalle2, ecalle3}, Hakim \cite{hakim1997transformations, hakim1998analytic},
Abate, Bracci, Tovena \cite{abate2001residual, ab_parabolic, abtann, ab}, Rong \cite{rong}, Molino \cite{molino},
Vivas \cite{vivas2011degenerate}, Arizzi, Raissy \cite{arra}, and others. A few generalizations to several variables of the
Leau-Fatou flower theorem have been proved, but none of them was strong enough to be able to describe the dynamics
in a full neighbourhood of the origin; furthermore, examples of unexpected phenomena not appearing in the one-dimensional
case have been found. Thus it is only natural to try and study the dynamics of meaningful classes of examples, with the aim
of extracting ideas applicable to a general setting; and Camacho's theorem suggests that a particularly interesting class
of examples is provided by time-1 maps of homogeneous vector fields. (Actually, the evidence collected so far strongly
suggests that a several variable version of Camacho's theorem might hold: generic germs tangent to the identity should be
locally topologically conjugated to time-1 maps of homogeneous vector fields. But we shall not pursue this topic here.)

This is the approach initiated in~\cite{ab}; and we discovered that there is a strong relationship between the dynamics
of the time-1 map of homogeneous vector fields and the dynamics of geodesics for meromorphic connections on Riemann surfaces.
To describe this relationship, we need to introduce a few notations and definitions.

Let $Q$ be a homogeneous vector field on $\C^n$ of degree~$\nu+1\ge 2$. First of all,
notice that the orbits of its time-1 map are contained in the real integral curves of~$Q$; so
we are interested in studying the dynamics of the \emph{real} integral curves of the \emph{complex}
homogeneous vector field~$Q$. (Actually, it turns out that complex integral curves of a homogeneous vector
field are related to --- classically studied--- sections which are horizontal with respect to a meromorphic
connection; see \cite{ab} for details.)

A \emph{characteristic direction} for $Q$
is a direction $v\in\PP^{n-1}(\C)$ such that the complex line (the \emph{characteristic leaf}) issuing from
the origin in the direction~$v$ is $Q$-invariant.  
An integral curve
issuing from a point of a characteristic leaf stays in that leaf forever; so the dynamics
in a characteristic leaf is one-dimensional, and thus completely known.
In particular,  if the vector field $Q$ is a multiple of the radial field (we shall say that $Q$ is \emph{dicritical})
every direction
is characteristic, the dynamics is one-dimensional and completely understood.
So, we are mainly interested in understanding the dynamics of integral curves outside the characteristic leaves of
non-dicritical vector fields.

Then in \cite{ab} we proved the following result:

\begin{theo}[Abate-Tovena \cite{ab}]
 Let $Q$ be a non-dicritical homogeneous
vector field of degree $\nu + 1 \geq 2$ in $\C^n$
and let $M_Q$ be the complement in $\C^n$ of the characteristic leaves of $Q$.
Let $[\cdot]\colon \C^n\setminus\{O\}\to\PP^{n-1}(\C)$ denote the canonical projection. Then there exists a
singular holomorphic foliation $\mathcal{F}$ of~$\PP^{n-1}(\C)$ in Riemann surfaces, and a partial meromorphic
connection $\nabla$ inducing a meromorphic connection on each leaf of~$\mathcal{F}$, whose poles coincide with the characteristic directions
of~$Q$, such that the following hold:
\begin{enumerate}
\item[(i)] if $\gamma\colon I \to M_Q$ is an integral curve of~$Q$ then the image of $[\gamma]$ is contained
in a leaf~$S$ of~$\mathcal{F}$ and it is a geodesic for~$\nabla$ in~$S$;
\item[(ii)] conversely, if $\sigma\colon I\to\PP^{n-1}(\C)$ is a geodesic for $\nabla$ in a leaf~$S$ of~$\mathcal{F}$
then there are exactly $\nu$ integral curves $\gamma_1,\ldots,\gamma_\nu\colon I\to
M_Q$ such that $[\gamma_j]=\sigma$ for $j=1,\ldots,\nu$.
\end{enumerate}
\end{theo}

Thanks to this result, we see that the study of integral curves for a homogeneous vector field in $\C^n$
is reduced to the study of geodesics for meromorphic connections on a Riemann surface $S$ (obtained
as a leaf of the foliation~$\Ff$).

In \cite{ab} the latter study has been carried out in the case $S=\PP^1(\C)$, which is the only case
arising when $n=2$ (and indeed it has led to a fairly extensive understanding of the dynamics of homogeneous
vector fields in $\C^2$, including the description of the dynamics in a full neighbourhood of the origin for
a substantial class of examples); the main goal of this paper is to extend this study from $\PP^1(\C)$ to
a generic compact Riemann surface, with the hope of applying in the future our results to the study of the
dynamics of homogeneous vector fields in $\C^n$ with $n\ge 3$. 

More precisely, we shall prove a Poincar\'e-Bendixson theorem (see Theorems~\ref{pbgeneral}
and~\ref{pbgeneral_res_grafi})
describing the $\omega$-limit set of the geodesics for a meromorphic connection
on a generic compact
Riemann surface. We recall that, in general, a ``Poincar\'e-Bendixson theorem"  is a result describing
recurrence properties of a class of
dynamical systems (see \cite{ciesielski2012} for a survey on the subject).
An important ingredient in our proof
will indeed be another Poincar\'e-Bendixson theorem, due to Hounie (\cite{hounie1981minimal}), describing the
minimal sets for smooth singular line fields on compact orientable surfaces (which in turn follows from a
similar statement for smooth vector fields, see \cite{pbgen}).

The paper is organized as follows.
In Section \ref{section_preliminary} we collect all the preliminary results that we shall need.
In Section \ref{section_meromorphic_connections} we generalize results obtained in \cite{ab} concerning geodesics
for meromorphic connections on the Riemann sphere to geodesics for meromorphic connections on a generic compact Riemann surface.
In Sections \ref{sezpbgeneral} and~\ref{resd} we prove our main Theorems:
we shall show that it is possible
to see a geodesic for a meromorphic connection
as part of an integral curve for a suitable line field on the surface,
which, in a neighbourhood of the support of the geodesic, is singular exactly on the poles of the connection;
then, thanks to
the result of Hounie mentioned above, we shall be able to give a topological description of
the possible minimal sets, and we shall conclude the proof by 
applying the theory developed in the previous sections. In Section~\ref{non-comp} we shall briefly
discuss what happens on non-compact Riemann surfaces.
Finally, in Section \ref{section_geod_toro} we present a detailed study of the geodesics
for holomorphic connections on complex tori (the only compact Riemann surfaces admitting meromorphic connections without poles).

\section{Preliminary notions}\label{section_preliminary}

Throughout this paper we shall denote by $S$ a Riemann surface
and by $p\colon TS\to S$ the tangent bundle on $S$. Furthermore, $\Oo_S$ will be the structure
sheaf of $S$, i.e., the sheaf of germs of holomorphic functions on $S$, and $\Mm_S$ the sheaf of germs of meromorphic functions.
$\Tt\Ss$ will denote the sheaf of germs of holomorphic sections of $TS$ and $\Mm_{TS}$ the sheaf of germs
of its meromorphic sections. Finally, $\Omega^1_S$ will denote the sheaf of germs of holomorphic 1-forms,
and $\Mm_S^1$ the sheaf of germs of meromorphic 1-forms.

We start recalling the definitions and the first properties of holomorphic and meromorphic connections
on $p\colon  TS \to S$. We refer to \cite{ab} and \cite{ilyako} for details.

\proc{Definition 2.1.}\label{defholconn}
A \emph{holomorphic} \emph{connection}
on the tangent bundle $p\colon TS\to S$ over a Riemann surface is a $\C$-linear
map $\nabla \colon \Tt\Ss \to \Omega^1_S \otimes \Tt\Ss$ satisfying
the Leibniz rule
\[
\nabla (fs) = df \otimes s + f \nabla s
\]
for all $s \in \Tt\Ss$ and $f \in \Oo_S$. We shall often write $\nabla_u s$ for $\nabla s(u)$.

A \emph{geodesic} for a holomorphic connection $\nabla$ is a real curve $\sigma \colon  I \to S$,
with $I \subseteq \R$ an interval, such that $\nabla_{\sigma'} \sigma' \equiv 0$.
\medbreak

Let us see what this definition means in local coordinates.
Given a holomorphic atlas $\{(U_\alpha,z_\alpha)\}$ on $S$, we denote by $\partial_\alpha:=\partial/\partial z_\alpha$ the induced local generator of $TS$ over $U_\alpha$. 
We shall always suppose that the $U_\alpha$'s
are simply connected.
On each~$U_\alpha$
we can find a holomorphic 1-form $\eta_{\alpha}$ such that
\[
\nabla \partial_\alpha = \eta_{\alpha} \otimes \partial_\alpha\;,
\]
and we shall say that the form $\eta_{\alpha}$ represents the
connection $\nabla$ on $U_\alpha$.
In fact we see that,
for a general section $s$ of $TS$, we can locally
compute $\nabla s$ by representing $s|_{U_\alpha}$ as
$s_\alpha \partial_\alpha$
for some holomorphic function $s_\alpha$ on $U_\alpha$ and writing
\[
\nabla (s_\alpha \partial_\alpha) = ds_\alpha \otimes \partial_\alpha + s_\alpha \nabla (\partial_\alpha) =
ds_\alpha \otimes \partial_\alpha + s_\alpha \eta_\alpha \otimes \partial_\alpha =
(ds_\alpha + \eta_\alpha) \otimes \partial_\alpha\;.
\]

\proc{Definition 2.2.}\label{defi_compatibile}
Let $g$ be an Hermitian metric on $p\colon  TS \to S$. We say that
a connection~$\nabla$ on $TS$ is \emph{adapted to}, or \emph{compatible with}, the metric $g$ if
\[
d \langle R,T \rangle= \langle \nabla R, T \rangle + \langle R,\overline{\nabla} T \rangle 
\]
for every pair of smooth vector fields $R$,~$T$,
that is if
\[
X\langle R,T \rangle = \langle \nabla_X R,T\rangle + \langle R, \nabla_{\bar{X}}T \rangle 
\]
for every triple of smooth vector fields $X$, $R$ and $T$ on $S$.
\medbreak

It is known that it is possible to associate to any Hermitian metric a
(not necessarily holomorphic) connection on the tangent bundle adapted to it, the
\emph{Chern connection}.
Locally, it is also possible to solve the converse problem, i.e., 
to construct local metrics adapted to a given holomorphic connection on $TS$.
To do this, it is convenient to consider
the local real function $n_\alpha \colon  U_\alpha \to \R^+$ given by
\[
n_\alpha (p)= g_p (\partial_\alpha, \partial_\alpha)\;.
\]
It is straightforward to see that $n_\alpha$ is a smooth function and, conversely,
we see that the function $n_\alpha$ uniquely determines the metric on $U_\alpha$.

With these notations in place it is not difficult to see that
the compatibility of a metric $g$ with a holomorphic connection $\nabla$ on
the domain $U_\alpha$ of a local chart is equivalent to 
\begin{equation}\label{dn=neta}
\partial n_\alpha = n_\alpha \eta_\alpha\;.
\end{equation}

Equation \eqref{dn=neta} can actually
be solved:
\begin{prop}[\cite{ab}, Proposition~1.1]\label{soluzden=neta}
Let $\nabla \colon  \Tt\Ss \to \Omega^1_S \otimes \Tt\Ss$
be a holomorphic connection on a Riemann surface $S$.
Let $\left( U_\alpha, z_\alpha\right)$ be a local chart
for~$S$ and let $\eta_\alpha \in \Omega^1_S (U_\alpha)$ be the 1-form representing $\nabla$ on $U_\alpha$. 
Given a holomorphic primitive $K_\alpha$ of $\eta_\alpha$ on $U_\alpha$, then
\[
n_\alpha = \exp (2\Re K_\alpha) = \exp (K_\alpha + \bar{K_\alpha})
\]
is a positive solution of \eqref{dn=neta}.

Conversely, if $n_\alpha$ is a positive solution of \eqref{dn=neta} then for any $z_0 \in U_\alpha$
and any simply connected neighbourhood $U\subseteq U_\alpha$ of $z_0$
there exists a holomorphic primitive $K_\alpha \in \Oo (U)$ of $\eta_\alpha$
over $U$ such that $n_\alpha = \exp (2 \Re K_\alpha)$ in U.
Furthermore, $K_\alpha$ is unique up to a purely imaginary additive constant.

Finally, two (local) solutions of \eqref{dn=neta} differ (locally) by a positive 
multiplicative constant.
\end{prop}

\proc{Remark 2.1.}
It is important to notice that Proposition~\ref{soluzden=neta} gives
only \emph{local} metrics adapted to~$\nabla$; a \emph{global} metric adapted to~$\nabla$ might not exist
(see \cite[Proposition~1.2]{ab}).
\medbreak

So we can associate to a holomorphic connection $\nabla$ a conformal family of compatible local metrics. It turns out (see \cite{ab})
that these local metrics are locally isometric to the Euclidean metric on $\C$.
In fact, given any holomorphic primitive $K_\alpha$ of $\eta_\alpha$, let the function
 $J_\alpha\colon  U_\alpha \to \C$ be
 a holomorphic primitive of $\exp({K_\alpha})$.
We immediately
remark that $J_\alpha$ actually exists, because $U_\alpha$ is simply connected, and
that
it is locally invertible, because $J_\alpha'= \exp ({K_\alpha})$.
In the following Proposition we summarize the main properties of $J_\alpha$.

\begin{prop}[\cite{ab}]\label{propJ}
Let $\nabla \colon  \Tt\Ss \to \Omega^1_S \otimes \Tt\Ss$
be a holomorphic connection on a Riemann surface $S$.
Let $\left( U_\alpha, z_\alpha\right)$ be a local chart
for $S$ with $U_\alpha$ simply connected, $\eta_\alpha$ the 1-form representing
$\nabla$ on~$U_\alpha$, and $K_\alpha$ a holomorphic primitive of $\eta_\alpha$
on~$U_\alpha$. Then every primitive $J_\alpha\colon U_\alpha\to\C$ of $\exp(K_\alpha)$ is a local 
isometry between $U_\alpha$, endowed with the metric represented by $n_\alpha=\exp(2\Re K_\alpha)$,
and $\C$, endowed with the Euclidean metric. In particular, every local metric adapted to~$\nabla$ is \emph{flat},
that is its Gaussian curvature vanishes identically. 

Moreover, a smooth curve $\sigma \colon I \to {U_\alpha}$ is a geodesic
for $\nabla$ if and only if there are two constants
$c_0$ and $w_0 \in \C$ such that 
$J_\alpha(\sigma(t) )= c_0 t + w_0$. In particular, the geodesic with
$\sigma(0)= z_0$ and $\sigma'(0)= v_0 \in \C^*$ is given by
\begin{equation}\label{eqJ1}
\sigma(t)= J_\alpha^{-1} \left(c_0 t + J_\alpha(z_0)\right)\;,
\end{equation}
where $c_0 = \exp \left( {K_\alpha}(z_0) \right) v_0$ and $J_\alpha^{-1}$ is the
analytic continuation of the local inverse of $J_\alpha$ near $J_\alpha(z_0)$
such that $J_\alpha^{-1} \left(J _\alpha(z_0)\right) = z_0$.

Finally, a curve $\sigma\colon [0,\eps)\to U_\alpha$ is a geodesic for $\nabla$ if and only if
\begin{equation}\label{eqJ2}
 \sigma' (t)= \exp\left( -{K_\alpha}(\sigma(t))\right) \exp\left({K_\alpha} (\sigma(0))\right)\sigma'(0)\;,
\end{equation}
if and only if
\begin{equation}\label{eqJ3}
 J_\alpha(\sigma(t)) = \exp\left( {K_\alpha} (\sigma(0))  \right) \sigma'(0) t + J_\alpha (\sigma(0)).
\end{equation}
\end{prop}
\medbreak

\proc{Remark 2.2.}\label{remangle}
Another consequence of the existence of a conformal family of local metrics adapted
to~$\nabla$ is that we can introduce on each tangent space $T_pS$ a well-defined
notion of \emph{angle} between tangent vectors, clearly independent of the particular local metric used to compute it.
\medbreak

We shall now give the official definition of \emph{meromorphic connection}.

\proc{Definition 2.3.}\label{defmeroconn}
A \emph{meromorphic connection} on the tangent bundle $p\colon TS\to S$ of a Riemann surface $S$
is a $\C$-linear map
$\nabla\colon \Mm_{TS} \to \Mm^1_S \otimes_{\Mm (S)} \Mm_{TS}$
satisfying the Leibniz rule
\[
\nabla (\tilde{f}\tilde{s}) = d\tilde{f} \otimes \tilde{s} + \tilde{f} \nabla \tilde{s}
\]
for all $\tilde{s} \in \Mm_{TS}$ and $\tilde{f} \in \Mm_S$.
\medbreak

It is easy to see that on a local chart $(U_\alpha,z_\alpha)$ a meromorphic
connection is represented by a meromorphic 1-form, that we continue to call $\eta_\alpha$, such that
\[
\nabla (\partial_\alpha ) = \eta_\alpha \otimes \partial_\alpha
\]
and so
\[
\nabla (\tilde s)= \nabla ( \tilde{s}_\alpha \partial_\alpha) = (d\tilde{s}_\alpha + \eta_\alpha )\otimes \partial_\alpha
\]
for every local meromorphic section $\tilde s|_{U_\alpha}= \tilde{s_\alpha} \partial_\alpha$ on $U_\alpha$,
exactly as it happens for holomorphic connections.

In particular, all the forms $\eta_\alpha$'s are
holomorphic if and only if $\nabla$ is a holomorphic connection.
More precisely, if we say that $p\in S$ is a \emph{pole} for a meromorphic
connection $\nabla$ if $p$ is a pole of~$\eta_\alpha$ for some (and hence any) local
chart $U_\alpha$ at~$p$, then $\nabla$ is a holomorphic connection on
the complement $S^0$ of the poles of $\nabla$ in~$S$.

This allows to define \emph{geodesics} for a meromorphic
connections as we did for holomorphic connections.

\proc{Definition 2.4.}
 A \emph{geodesic} for a meromorphic connection $\nabla$ on $p\colon TS \to S$
 is a real curve $\sigma\colon  I \to S^0$, with $I \subseteq \R$ an interval, such that $\nabla_{\sigma'} \sigma' \equiv 0$.
\medbreak

The following definition will be of primary importance in the sequel.

\proc{Definition 2.5.}
The \emph{residue} $\Res_p\nabla$ of a meromorphic connection $\nabla$ at
a point $p \in S$ is the residue of any 1-form $\eta_\alpha$ representing $\nabla$ on a local chart $(U_\alpha,z_\alpha)$ at~$p$.
Clearly, $\Res_p\nabla\ne 0$ only if $p$ is a pole of~$\nabla$.
\medbreak

It is easy to see
that the residue of a meromorphic connection at a point $p \in S$
is well defined, i.e., it does not depend on the particular chart~$U_\alpha$. 
Moreover, by the Riemann-Roch formula, the sum of the residues of a meromorphic connection $\nabla$ is independent from
the particular connection, as stated in the following classical Theorem (see \cite{ilyako} for a proof).

\begin{theo}\label{sum=-chi}
Let $S$ be a compact Riemann surface and $\nabla$ a meromorphic
connection on $p\colon  TS \to S$. Then
\begin{equation}\label{eqsum=chi}
\sum_{p \in S} \Res_p \nabla = - \chi_S\;,
\end{equation}
where $\chi_S$ is the Euler characteristic of $S$.
\end{theo}

Our main result is a classification of the possible $\omega$-limit sets
for the geodesics of meromorphic connections on the tangent bundle of a compact Riemann surface.
We recall that the $\omega$-limit set $\omega (\sigma)$ of a curve $\sigma\colon  [0,\eps)\to S^0$ is given by
the points $p$ of $S$ for which there exists a sequence $\{t_n\}$, with $t_n \to \eps$,
such that $\sigma(t_n)\to p$. It is easy to see that
\[
\omega (\sigma) = \bigcap_{\eps' < \eps} \overline{ \{ \sigma(t)\colon t > \eps'\}}.
\]

The main part of the proof will consist in extending the
tangent field of the geodesic to a smooth line field (a rank-1 real foliation) on $S$, singular only on
the poles of the connection.
Then we shall be able to use the following theorem by Hounie (\cite{hounie1981minimal}) to study the minimal
sets of this line field, where a \emph{minimal set} for a line field $\Lambda$ is
a closed, non-empty, $\Lambda$-invariant subset of $S$ 
without proper subsets having the same properties:

\begin{theo}[Hounie \cite{hounie1981minimal}]\label{pbfoliaz}
Let $S$ be a compact connected two-dimensional smooth real manifold (e.g., a
Riemann surface)
and let $\Lambda$ be a smooth line field with singularities on $S$. Then
a $\Lambda$-minimal set $\Omega$ must be one of the following:
 \begin{enumerate}
  \item a singularity of $\Lambda$;
  \item a closed integral curve of $\Lambda$, homeomorphic to $S^1$;
  \item all of $S$, and in this case $\Lambda$ is equivalent to an irrational line field
  on the torus (i.e., there exists a homeomorphism
  $\phi \colon  S \to T$ between $S$ and a torus $T$ transforming the given foliation into one induced by
  an irrational line field).
 \end{enumerate}
\end{theo}

\section{Meromorphic connections on the tangent bundle}\label{section_meromorphic_connections}

In this section, we study in more detail geodesics for meromorphic connections on the tangent bundle of a compact Riemann surface.
In particular,
we extend results contained in Section 4 of \cite{ab} from $\PP^1(\C)$ to the case of a generic compact Riemann surface $S$.

To do so, we start introducing the following definitions/notations.

\proc{Definition 3.1.}\label{defi_cicli_geod}
Let $S$ be a compact Riemann surface. Let $\nabla$ be a meromorphic connection
on $TS$ and let $S^0\subseteq S$ be the complement
of the poles.
\begin{itemize}
\item A \emph{geodesic ($n$-)cycle} is the union of
$n$ geodesic segments $\sigma_j\colon [0,1]\to S^0$, disjoint except for the conditions
$\sigma_j (0) = \sigma_{j-1} (1)$ for $j=2,\ldots,n$ and $\sigma_{1} (0)= \sigma_n (1)$. The points $\sigma_j(0)$ will be
called \emph{vertices} of the geodesic cycle. 
\item A \emph{($m$-)multicurve} is a union of $m$ disjoint geodesic cycles.
A multicurve will be said to be \emph{disconnecting} if it disconnects $S$, \emph{non-disconnecting} otherwise.
\item A \emph{part} $P$ is the closure of 
a connected open subset of $S$ whose boundary is a multicurve~$\gamma$.
A component $\sigma$ of $\gamma$ is \emph{surrounded} if the interior of $P$
contains both sides of a tubular neighbourhood of~$\sigma$ in~$S$; it is \emph{free} otherwise.  
The \emph{filling} $\hat P$ of a part $P$ is the compact surface obtained
by glueing a disk along each of the free components of~$\gamma$, and  
not removing any of the surrounded components of~$\gamma$. 
\end{itemize}
\medbreak

We remark that a part may be all of $S$, when the associated multicurve is
non-disconnecting (which is equivalent to saying that the multicurve has no free components).
Moreover, we see that every disconnecting multicurve contains the boundary
of a part $P \varsubsetneq S$.

As recalled in the previous section, we can associate to a meromorphic
connection~$\nabla$ conformal families of local metrics on the regular part
$S^0\subseteq S$ of the Riemann surface, and we have a well-defined notion of
angle between tangent vectors. In particular, it makes sense to speak
of the \emph{external angle}~$\epsilon$ at a vertex $\sigma_j(0)$ of a geodesic cycle as the angle between the tangent vectors 
$\sigma'_{j-1}(1)$ and $\sigma'_j(0)$; by definition, $\epsilon\in(-\pi,\pi)$.

Again using the conformal families of local metric it is possible to define the \emph{geodesic curvature}
$k_g$ of a curve contained in~$S^0$. 
In particular, if $p\in S$ is a pole of~$\nabla$ and $\tau$ is a small (clockwise) circle around $p$, not containing
other poles, we have
\begin{equation}\label{kg1polo}
\int_{\tau} k_g = -2\pi (1 +  \Re \Res_p (\nabla) )\;;
\end{equation}
see \cite[Theorem 4.1]{ab}.

With all these ingredients at our disposal, it becomes natural to try
and apply a Gauss-Bonnet Theorem to study the relation between the residues
of the connection~$\nabla$ in a part $P$ of $S$,
the external angles
at the vertices of the multicurve bounding $P$ and the topology of $S$ (cp. \cite[Theorem 4.1]{ab}).

\begin{theo}\label{ab4.1g}
Let $\nabla$ be a meromorphic connection on a compact Riemann surface~$S$, with poles
$\{p_1, \dots ,p_r \}$ and set $S^0 :=S \setminus \{p_1, \dots ,p_r \}$.
Let $P$ be a part of $S$ whose boundary multicurve $\gamma\subset S^0$ has
$m_f\ge 1$ free components, positively oriented with respect to~$P$. 
Let $z_1,\dots , z_s$
denote the vertices of the free components of~$\gamma$, and $\eps_j\in(-\pi,\pi)$ the
external angle at $z_j$. Suppose that $P$ contains the poles
$\{p_1, \dots ,p_g \}$ and denote by $g_{\hat P}$ the genus of the filling $\hat P$ of~$P$.
Then
\begin{equation}\label{eqGB}
\sum_{j=1}^{s} \eps_j = 2 \pi \left( 2 - m_f -  2g_{\hat P} + \sum_{j=1}^{g} \Re \Res_{p_j}(\nabla) \right)\;.
\end{equation}
\end{theo}

\proc{Proof.}
For $j=1,\ldots,g$ let $\tau_j$ be a small clockwise circle bounding a disk in~$P$ centered at~$p_j$,
and let $k_g^j$ be the geodesic curvature of $\tau_j$.

Applying the Gauss-Bonnet Theorem
as in \cite[Theorem 4.1]{ab}
to the complement $P^0$ in $P$ of the small disks bounded by $\tau_1,\ldots,\tau_g$ we find that
\begin{equation}\label{eq1_4.1}
\sum_{j=1}^{g} \int_{\tau_j} k_g^j + \sum_{j=1}^{s} \eps_j = 2 \pi \chi_{P^0} = 2 \pi (2 - m_f- g - 2g_{\hat P})\;.
\end{equation}
But from \eqref{kg1polo} we get
\begin{equation}\label{eq2_4.1}
\sum_{j=1}^{g} \int_{\tau_j} k_g^j= -2\pi g - 2\pi \sum_{j=1}^{g} \Re \Res_{p_j}(\nabla)\;.
\end{equation}
Comparing \eqref{eq1_4.1} and \eqref{eq2_4.1} we get the assertion.
\ep
\medbreak

\proc{Remark 3.1.}
When the multicurve $\gamma$ does not disconnect $S$ then Theorem~\ref{ab4.1g}
reduces to Theorem~\ref{sum=-chi}. Indeed, in this case $\gamma$ has no free components, and thus \eqref{eqGB} becomes
\[
0= 2 \pi \left( 2 -  2g_{S} + \sum_{j=1}^{g} \Re \Res_{p_j}(\nabla) \right)\;,
\]
which is equivalent to \eqref{eqsum=chi}.
\medbreak

In the next two Corollaries we highlight what happens when the disconnecting multicurve is made up by a single geodesic or by
a single geodesic cycle composed by two geodesics (cp. \cite[Corollaries 4.2 and 4.3]{ab}).

\begin{cor}[One disconnecting geodesic]\label{ab4.2g}
Let $\nabla$ be a meromorphic connection on a compact Riemann surface $S$, with poles
$\{p_1, \dots ,p_r \}$ and set $S^0:=S \setminus \{p_1, \dots ,p_r \}$.
Let $\sigma$ be a
disconnecting geodesic $1$-cycle.
Let $P$ be one of the two parts in which $S$ is disconnected
by $\sigma$ and $\eps \in (-\pi, \pi)$ the unique external angle of $\sigma$. Then
\[
\eps = 2 \pi \left( 1 - 2 g_{\hat P} + \sum_{p_j\in P} \Re \Res_{p_j}(\nabla) \right)\;.
\]
In particular,
\[
\sum_{p_j \in P} \Re \Res_{p_j}(\nabla) \in (-3/2 + 2 g_{\hat P}, -1/2 + 2 g_{\hat P})\;.
\]
\end{cor}

\begin{cor}[Two geodesics whose union disconnects $S$]\label{ab4.3g}
Let $\nabla$ be a meromorphic connection on a compact Riemann surface $S$, with poles
$\{p_1, \dots ,p_r \}$ and set $S^0:=S \setminus \{p_1, \dots ,p_r \}$.
Let $\gamma$ be a disconnecting geodesic $2$-cycle.
Let $P$ be one of the two parts in which $S$ is disconnected
by $\gamma$, and $\eps_0$, ~$\eps_1\in(-\pi,\pi)$ the two external angles of $\gamma$. Then
\[
\eps_0 + \eps_1= 2 \pi \left( 1 - 2 g_{\hat P} + \sum_{p_j\in P} \Re \Res_{p_j}(\nabla) \right)\;,
\]
and hence
\[
\sum_{p_j \in P} \Re \Res_{p_j}(\nabla) \in (-2 + 2g_{\hat P} , 2 g_{\hat P})\;.
\]
\end{cor}

\proc{Remark 3.2.}\label{rempo}
In particular, a part of $S$ bounded by a disconnecting $1$-cycle
or by a disconnecting $2$-cycle must necessarily contain a pole,
because $0\notin (-2 + 2g_{\hat P} , 2 g_{\hat P})$.
\medbreak

In the following we shall need to consider \emph{closed} geodesics and \emph{periodic} geodesics for a meromorphic connection.

\proc{Definition 3.2.}
 A geodesic $\sigma\colon  [0,l] \to S$ is \emph{closed} if $\sigma(l) =\sigma(0)$ and $\sigma' (l)$ is a positive multiple
 of $\sigma' (0)$; it is \emph{periodic} if $\sigma(l)=\sigma(0)$ and $\sigma' (l)= \sigma'(0)$.
\medbreak

By Corollary \ref{ab4.2g} we immediately see that a disconnecting
geodesic 1-cycle $\sigma$ is a closed geodesic if and only if for every part $P$ of $S$ bounded by $\sigma$ we have
\begin{equation}\label{formula_closed_generica}
\sum_{p_j\in P} \Re \Res_{p_j}(\nabla) = -1 + 2 g_{\hat P}\;.
\end{equation}

Contrarily to the Riemannian case, closed geodesics are not necessarily
periodic; for examples (and more) in the Riemann sphere see \cite{ab}, and for
examples on a complex torus see Section \ref{section_geod_toro}. 

\section{$\omega$-limits sets of geodesics}\label{sezpbgeneral}

In this section we prove our main Theorem, the classification of the possible $\omega$-limit sets of geodesics
for a meromorphic connection on a compact Riemann surface.
The main idea will be to see a not self-intersecting geodesic $\sigma$ as part of an integral curve of a
suitable line field $\Lambda$ on $S$, and then to apply
Theorem \ref{pbfoliaz} to get some information about the minimal
sets for $\Lambda$
contained in the $\omega$-limit set
of $\sigma$. Then we shall use these information to discuss the shape of the $\omega$-limit set itself.

Let us start by studying the local structure of the $\omega$-limit set of a not self-intersecting geodesic~$\sigma$.

\begin{prop}\label{local_omega}
 Let $S$ be a Riemann surface,
 $\nabla$ a meromorphic connection on~$S$, and $S^0$
 the complement of the poles for $\nabla$.
 Let $\sigma \colon I \to S^0$ be a maximal not self-intersecting geodesic
 for $\nabla$ and denote by $W$ its $\omega$-limit set.
 Let $z \in W\cap S^0$. Then:
 \begin{enumerate}
  \item[(i)] there exists a unique direction $v$ at $z$ such that
  for every sequence $\{z_n\} \subset \sigma(I)$ with $z_n\to z$
the directions $\sigma'(z_n)$ converge to $v$; moreover, the image
of the maximal geodesic~$\sigma_z$ issuing from~$z$ in the direction~$v$ is contained in~$W$;
\item[(ii)] if $W$ contains a curve $\tau$ passing through~$z$
transversal to~$\sigma_z$ then~$z \in \mathring W$;
\item[(iii)] if $z\in \de W$ then $\sigma_z$ is the unique
geodesic segment through $z$ contained in~$W$.
 \end{enumerate}
\end{prop}

\proc{Proof.}
(i) The existence of $v$ follows from the fact that $\sigma$ does not
self-intersect. Using a local isometry $J$ defined in a neighbourhood of~$z$
transforming geodesic segments in Euclidean segments we immediately see
that the image of the geodesic segment issuing from $z$ in the direction~$v$
must be contained in~$W$, and the maximality follows from the maximality of~$\sigma$.
\smallskip

(ii) Up to using a local isometry~$J$ we can assume that all geodesic
segments in a neighbourhood of~$z$ are Euclidean segments. 

Consider a point $z'\ne z$ belonging to~$\tau$. Since $z' \in W \subset S^0$, (i)
gives us
a geodesic
segment  $\sigma_{z'}$ through $z'$ contained in $W$. Notice that $\sigma_{z'} \neq \tau$
as soon as~$z'$ is close enough to~$z$. In fact,
since we have segments of $\sigma$ accumulating $\sigma_z$, the
segments of~$\sigma$ accumulating $\sigma_{z'}$ cannot converge to $\tau$ without
forcing $\sigma$ to self-intersect, impossible. The mapping $z'\mapsto\sigma_{z'}$
is continuous, again because $\sigma$ cannot self-intersect; thus the geodesic
segments $\sigma_{z'}$ fill an open neighbourhood of~$z$ contained in~$W$, and $z\in\mathring W$ as claimed.
\smallskip

(iii) It immediately follows from (ii).
\ep
\medbreak

\proc{Definition 4.2.}
Let $S$ be a Riemann surface, $\nabla$ a meromorphic connection on~$S$, and $S^0$
 the complement of the poles of $\nabla$.
 Let $\sigma \colon I \to S^0$ be a maximal not self-intersecting geodesic
 for $\nabla$, with $\omega$-limit~$W\subseteq S$. Given $z\in W\cap S^0$,
 the maximal geodesic $\sigma_z$ given by Proposition~\ref{local_omega}.(i) issuing
 from $z$ and contained in~$W$ is the \emph{distinguished geodesic} issuing from~$z$.
\medbreak

The following two lemmas contain basic properties of distinguished geodesics.

\begin{lemma}\label{nointersect}
Let $S$ be a Riemann surface, $\nabla$ a meromorphic connection on~$S$, and $S^0$
 the complement of the poles of $\nabla$.
 Let $\sigma \colon I \to S^0$ be a maximal not self-intersecting geodesic
 for $\nabla$, with $\omega$-limit~$W\subseteq S$. Let $\sigma_z$ be the
 distinguished geodesic issuing from $z\in W\cap S^0$. Then either
 $\sigma_z=\sigma$ or $\sigma$ does not intersect~$\sigma_z$.
\end{lemma}

\proc{Proof.}
If $\sigma$ intersects $\sigma_z$, either $\sigma=\sigma_z$ or the intersection is transversal;
but $\sigma$ accumulates~$\sigma_z$, and thus in the latter case $\sigma$ must self-intersect,
contradiction. 
\ep
\medbreak

\begin{lemma}\label{distclosed}
Let $S$ be a Riemann surface, $\nabla$ a meromorphic connection on~$S$, and $S^0$
 the complement of the poles of $\nabla$.
 Let $\sigma \colon I \to S^0$ be a maximal not self-intersecting geodesic for $\nabla$, with $\omega$-limit~$W\subseteq S$.
 Assume there exists a closed distinguished geodesic~$\alpha\subseteq W$.
 Then $W=\alpha$. 
 \end{lemma}
 
\proc{Proof.} If $\alpha=\sigma$ then the assertion is trivial; so
by the previous lemma we can assume that $\sigma$ does
not intersect~$\alpha$.

For each $w\in\alpha$ we can find a neighbourhood $A_w\subset S^0$ of~$w$,
contained in a tubular neighbourhood of~$\alpha$, and an isometry $J_w\colon A_w\to B_w$ such that:
\begin{itemize}
\item[(a)] $B_w$ is a quadrilateral containing the origin, with two opposite sides parallel to the imaginary axis; 
\item[(b)] $J_w(w)=O$ and $J_w(\alpha\cap A_w)$ is the intersection of~$B_w$ with the real axis;
\item[(c)] the $J_w$-images of the connected components of the intersection
of~$\sigma$ with~$A_w$ contained in the lower half-plane intersect both
vertical sides of~$B_w$ and accumulate the real axis (this can be achieved using Proposition~\ref{local_omega}.(i)).
\end{itemize}
By compactness, we can cover $\alpha$ with finitely many $A_{w_0},\ldots,A_{w_r}$
of such neighbourhoods; furthermore, we can also assume that
each $A_{w_j}$ intersects only $A_{w_{j-1}}$ and $A_{w_{j+1}}$, with the
usual convention $w_{-1}=w_r$ and $w_{r+1}=w_0$;
let $A=A_0\cup\cdots\cup A_r$. Furthermore, since $A$ is contained in a
tubular neighbourhood of~$\alpha$ and $S$ is orientable, $A\setminus\alpha$ has
exactly two connected components; and we can assume that the segments of $\sigma$
mapped by the $J_{w_j}$ in the lower half-plane (and thus accumulating~$\alpha$)
are all contained in the same connected component of~$A\setminus\alpha$.

Let $\tau$ be the segment of geodesic issuing from~$w_0$ such that $J_{w_0}(\tau)$ is the
intersection of~$B_{w_0}$ with the negative imaginary axis. Denote by $\{z_j\}$ the infinitely many intersections
of $\sigma$ with $\tau$, accumulating~$w_0$. Again by compactness, we can assume that if we follow~$\sigma$ starting from a~$z_1$
close enough to~$w_0$ we stay in~$A$ until we get to the next intersection~$z_2$.
Furthermore, by property (c) the segment $\sigma_1$ of~$\sigma$ from $z_1$ to $z_2$ must
intersects all $A_{w_k}$'s either in clockwise or in counter-clockwise order. Let $D\subset A$ denote
the domain bounded by $\alpha$, the geodesic segment $\sigma_1$ and the geodesic segment of $\tau$ from $z_1$ to~$z_2$.  

Assume first that $z_2$ is closer to $w_0$ than $z_1$. In this case, if we follow $\sigma$ starting
from $z_2$ property (c) forces $\sigma$ to remain inside $D$ because it cannot
intersect itself nor $\alpha$. Furthermore, property (c) again implies that the
next intersection $z_3$ should be closer to~$w_0$ than~$z_2$. We can then repeat
the argument, and we find that successive intersections of~$\sigma$ with~$\tau$ monotonically converge to~$w_0$. This implies that 
$\sigma$ accumulates $\alpha$ and nothing else, that is $W=\alpha$ as claimed.

If instead $z_2$ is farther away from $w_0$ than $z_1$, if we follow $\sigma$ starting
from $z_2$ we may possibly leave~$A$; but since $\sigma$ accumulates~$w_0$ we must sooner or
later get back to~$D$, and the only way is intersecting $\tau$ in a point $z_3$ between $z_1$ and $z_2$. But
now following $\sigma$ starting from $z_3$ we are forced to stay into~$D$ until we intersect $\tau$ in a point $z_4$ necessarily
closer
to~$w_0$ than $z_1$ --- and a fortiori closer than $z_3$.
We can then repeat the previous argument using $z_3$ and $z_4$ instead of~$z_1$ and $z_2$ to get again the assertion.
\ep
\medbreak

The main tool for the study of the global behavior of not self-intersecting
geo\-desics is the following Theorem, providing a smooth line field
$\Lambda$ such that $\sigma$ is (part of) an integral curve of $\Lambda$.

\begin{theo}\label{pblemma}
Let $S$ be a Riemann surface
and $\nabla$ a meromorphic connection on $S$,
with poles $p_1, \dots, p_r \in S$. Let $S^0=S \setminus \{ p_1, \dots, p_r \}$.
Let $\sigma \colon I \to S^0$ be a geodesic for $\nabla$ without
self-intersections, maximal in both forward and backward time. 
Then there exists a smooth line field $\Lambda$ with singularities on $S$ which has
$\sigma$ as integral curve and, in a
neighbourhood of $\sigma$, is singular exactly on the poles of $\nabla$.
Furthermore, the $\omega$-limit~$W$ of $\sigma$ is 
$\Lambda$-invariant, and $\Lambda|_W$ is uniquely determined.
\end{theo}

\proc{Proof.}
If $A \subset S^0$ is open, simply connected and small enough, we can find a metric~$g$
on $A$ compatible with $\nabla$ and an isometry $J$ between $A$ endowed with $g$
and an open set in $\C$ endowed with the euclidean metric; $g$ is unique
up to a positive multiple (see Propositions \ref{soluzden=neta} and \ref{propJ}).

We consider an open cover $\mathcal{A}=\{A_i\}$ of $S^0$ with
the following properties:
\begin{enumerate}
 \item[(a)] $\Aa$ is locally finite;
 \item[(b)] each $A_i$ is endowed with a metric $g_i$ compatible with $\nabla$, and with an isometry
 $J_i \colon  A_i \to B_i \subseteq \C$, with $B_i$ convex;
 \item[(c)] if $\sigma$ intersects $A_i$ then the angle between any pair of
 lines in $\C \supset B_i$ containing a component (which necessarily is a line segment) of the $J_i$-image of $\sigma$
is less than $\frac{\pi}{8}$. This can be achieved thanks to the smooth dependence of the
solution of the geodesic equation from the initial conditions and the fact that the isometry is smooth.
\end{enumerate}

We shall start building a line field on every open set $A_i$ of $\Aa$. Then we shall show how to use them to find
a global line field $\Lambda$ on $S^0$ having~$\sigma$ as integral curve, and finally we shall extend $\Lambda$ to all of $S$.

Let us then choose an open set $A \in \Aa$, together with its image $B$.
We shall now construct a smooth flow of curves in $B$ that
will correspond to a smooth flow of curves, and so
to a line field, in $A$.

If $\sigma$ does not cross $A$, we put on $B$ any smooth vector
field which is never zero and consider its associated (regular) foliation.

If $\sigma$ crosses $A$, we consider the segments
$\sigma_n \subset B$ which are the images of the connected components of
the intersection of~$\sigma$ with~$A$ (recall that $J$ sends geodesics segments in $A$ to Euclidean segments in $B$).
By our assumption on $A$, the angle between
$\sigma_i$ and $\sigma_j$ is bounded by $\frac{\pi}{8}$ for every pair $(i,j)$.

By the convexity of $B$ and the maximality of~$\sigma$, the $\sigma_n$'s subdivide $B$ in connected components. We describe now
how to costruct the flow in all these components.
\begin{enumerate}
\item\label{blato} If a connected component $C$ is open and its boundary
(in $B$) consists of a unique~$\sigma_0$, then we define our flow on $C$ by
means of lines parallel to $\sigma_0$ and take the associated line field.
\item\label{b01} If a connected component $C$ is open and it is bounded by
two segments $\sigma_0$ and~$\sigma_1$ we define the vector field in the following way:
we take two points $y_0 \in \sigma_0$ and $y_1 \in \sigma_1$. By convexity of $C$, the segment joining them
is contained in~$\bar{C}$. We parametrize this segment as $\tau \colon [0,1] \to\bar{C}$ with $\tau(0)=y_0$ and $\tau(1)=1$. 
For every $t \in [0,1]$, we consider the line $l_t$ passing through~$\tau(t)$ forming an angle
$\phi(t)\theta_1$ with $\sigma_0$, where $\phi \colon [0,1] \to [0,1]$ is a suitable smooth not decreasing
function which is 0 in a
neighbourhood of 0 and $1$ in a neighbourhood of $1$, and $\theta_1$ is the angle between $\sigma_0$ and~$\sigma_1$.
We immediately
see that the intersections $\tilde{l}_t = l_t \cap B$ form a smooth
flow on $C$, and that this flow is smooth also at the boundary of $C$ (i.e.\, near $\sigma_0$ and $\sigma_1$).
\item\label{limitflow}
Since the $\sigma_n$'s are disjoint, maximal and with angles bounded by $\pi/8$, 
the only missing case is a component $C$ consisting of a segment $\tau$
accumulated by a subsequence of $\sigma_n$; we then add this segment to the line field. \end{enumerate}
\smallbreak

In this way, thanks to the smooth dependence of geodesics on initial conditions,
we have defined a smooth and never vanishing line field on $B$ and hence,
via $J$, on~$A$. Notice that the image in~$B$ of a component of the
intersection of the $\omega$-limit set~$W$ with~$A$ must be a segment as
in case~\ref{limitflow}; thus $W\cap A$ is invariant with respect to this local line field,
which is uniquely defined there.

The next step consists in glueing the local line fields we have built on the $A_i$ to a global field on $S^0$.
This means that we must specify,
for every point $p \in S^0$, a direction $\lambda (p)$ in $T_pS^0$ such that the correspondence
$p \mapsto \lambda (p)$ is smooth.
To do so, we consider a partition of unity $\{\rho_i\}$ subordinated to the cover $\Aa$. If $p$ belongs to a
unique $A_i$, we use as $\lambda(p)$ the one given by the local costruction
above. Otherwise, if $p$ belongs to a finite number
of $A_i$'s (recall that the cover is locally finite) we do the following.
Suppose that $p \in A_1 \cap \dots \cap A_n$. We have $n$ lines in $T_pS^0$,
given by the local constructions on the $B_i$'s.
We use the partition of unity to do a convex combination of (the angles with
respect to any fixed line of) these lines, thus obtaining a line in $T_pS^0$.
We remark that, by Remark~2.2, the angles are the same on all $A_1,\ldots,A_n$,
and the bounds on the angles ensure that the convex combination yields a well-defined line.

We have thus obtained a smooth line field on $S^0$ having $\sigma$ as integral
curve. Since it is non-singular along~$\sigma$, by smoothness it is non-singular
in a neighbourhood of~$\sigma$ in~$S^0$. 

Finally, we extend this line field to all of $S$, adding the poles as singular points.
The resulting field satisfies the requests of the Theorem,
and we are done.
\ep
\medbreak

\proc{Remark 4.1.}\label{rem41}
The reason we have to use a line field instead of a vector field is step~3 in the
previous proof: a priori, the geodesic $\sigma$ might accumulate the segment~$C$
from both sides along opposite directions, preventing the identification of a non-singular vector field along~$C$.
\medbreak

\proc{Definition 4.1.}\label{def41}
Let $S$ be a Riemann surface,
$\nabla$ a meromorphic connection on~$S$, and
$S^0\subset S$ the complement of the poles.
Let $\sigma \colon I \to S^0$ be a maximal geodesic for $\nabla$ without
self-intersections, maximal in both forward and backward time. 
A smooth line field $\Lambda$ with singularities on $S$ as in the statement of
Theorem~\ref{pblemma} is called \emph{associated} to~$\sigma$.
\medbreak

The uniqueness of an associated line field on the $\omega$-limit set suggests
the following definition. We recall that a \emph{minimal set} for a line
field $\Lambda$ is a minimal element (with respect to inclusion) in the family of closed $\Lambda$-invariant subsets.

\proc{Definition 4.3.}\label{def42}
Let $S$ be a Riemann surface, $\nabla$ a meromorphic connection on~$S$, and $S^0$
 the complement of the poles of $\nabla$.
 Let $\sigma \colon I \to S^0$ be a maximal not self-intersecting geodesic
 for $\nabla$, with $\omega$-limit set~$W\subseteq S$. Then 
a \emph{minimal set} for $\sigma$ is a minimal set contained in~$W$ of any line field associated to~$\sigma$.
\medbreak

The following result characterizes the possible minimal sets for a maximal not
self-intersecting geodesic in a compact Riemann surface.

\begin{theo}\label{prop_minimali_limite}
 Let $S$ be a compact
 Riemann surface and $\nabla$ a meromorphic connection on $S$. Let $S^0$
 be the complement of the poles of $\nabla$. Let $\sigma \colon [0, \eps_0) \to S^0$ be
 a maximal not self-intersecting geodesic for $\nabla$.
 Then the possible minimal sets for $\sigma$ are the following:
 \begin{enumerate}
  \item a pole of $\nabla$;
  \item\label{minimale_sigma_chiusa} a closed curve, homeomorphic to $S^1$, which
  is a closed geodesic for $\nabla$, and in this case the $\omega$-limit set of
  $\sigma$ coincides with this closed geodesic;
  \item all of $S$, and in this case $S$ is a torus and $\nabla$ is holomorphic everywhere.
 \end{enumerate}
\end{theo}

\proc{Proof.}
Call $p_0=\sigma(0)$ the starting point of $ \sigma$, and $W$ the $\omega$-limit set of~$\sigma$.

We apply the construction in Theorem~\ref{pblemma}, considering $p_0$ as a virtual
pole. In this way, $\sigma$ becomes maximal in both
forward and backward time and the construction can be carried out as before.
We thus build a line field $\Lambda$ with singularities on $S$, which in a neighbourhood of $\sigma$
is singular
exactly on the poles of $\nabla$ and on $p_0$, and having $\sigma$ as integral curve.

Applying Theorem~\ref{pbfoliaz} to~$\Lambda$ we see that the minimal sets for $\sigma$ can be: all
of~$S$ (and in this case $S$ is a torus and $\nabla$ has no poles), a closed curve
homeomorphic to $S^1$ (which is a closed geodesic thanks to
Proposition \ref{local_omega}), a pole of $\nabla$ or $p_0$. This last possibility
is excluded by considering a new starting point $p_0'$ for $\sigma$ close enough
to~$p_0$ and noticing that now $p_0$ is not a pole for the new line field.

Finally, the closed geodesic in case~2 being $\Lambda$-invariant is necessarily
distinguished; we can then end the proof by quoting Lemma~\ref{distclosed}.
\ep
\medbreak

To describe the possible $\omega$-limit sets of a $\nabla$-geodesic, we need
a definition.

\proc{Definition 4.3.}
A \emph{saddle connection} for a meromorphic connection on a
Riemann surface $S$ is a maximal geodesic $\sigma\colon (-\eps_{-}, \eps_{+}) \to S^0$ such that $\sigma (t)$ tends to a pole for
both $t \to -\eps_-$ and $t\to\eps_+$.

A \emph{graph of saddle connections} is a connected planar
graph in $S$ whose vertices are poles and whose arcs are disjoint saddle connections. A \emph{spike} is a saddle connection
of a graph which does not belong to any cycle of the graph.

A \emph{boundary graph of saddle connections} (or \emph{boundary graph}) is a graph of saddle connections
which is also the boundary of a connected open
subset of $S$. A boundary graph is \emph{disconnecting} if its complement in $S$ is not connected.
\medbreak

We can now state and prove our main result.

\begin{theo}\label{pbgeneral}
Let $S$ be a compact Riemann surface and $\nabla$ a meromorphic connection on $S$,
with poles $p_1, \dots, p_r \in S$. Let $S^0=S \setminus \{ p_1, \dots, p_r \}$.
Let $\sigma \colon [0, \eps_0) \to S^0$ be a maximal geodesic for $\nabla$.
Then either
\begin{enumerate}
\item\label{pbtendpole} $\sigma(t)$ tends to a pole of $\nabla$ as $t\to \eps_0$; or
\item\label{pbclosed} $\sigma$ is closed; or
\item\label{pbtendclosed} the $\omega$-limit set of $\sigma$ is the support of a closed geodesic; or
\item\label{pbtendcycle} the $\omega$-limit set of $\sigma$ in $S$ is a boundary graph of saddle connections; or
\item\label{pbtutto} the $\omega$-limit set of $\sigma$ is all of $S$;
\item\label{pbinterno} the $\omega$-limit set of $\sigma$ has non-empty interior and non-empty boundary, and each component of
its boundary is a graph
of saddle connections with no spikes and at least one pole (i.e., it cannot be a closed geodesic); or
\item\label{pbinfint} $\sigma$ intersects itself infinitely many times.
\end{enumerate}
Furthermore, in cases \ref{pbclosed} or \ref{pbtendclosed} 
the support of $\sigma$ is contained in only one of the components
of the complement of the $\omega$-limit set, which is a part $P$ of $S$ having the $\omega$-limit set as boundary. 
\end{theo}

\proc{Proof.}
Suppose $\sigma$ is not closed, nor with infinitely many self-intersections.
Then up to
changing the starting point of $\sigma$
we can assume that $\sigma$ does not self-intersect. Let $W$ be the $\omega$-limit set of~$\sigma$.

By Zorn's lemma, $W$ must contain a minimal set, that 
by Proposition \ref{prop_minimali_limite} must be either 
a pole, a closed geodesic
or all of $S$ (and in this case $W=S$, which gives~\ref{pbtutto}).

If $W$ reduces to a pole or to a closed geodesic,
we are in cases \ref{pbtendpole} or \ref{pbtendclosed}.
So we can assume that $W$ properly contains a minimal set and is properly contained in~$S$; we have
to prove that we are in case \ref{pbtendcycle}
or \ref{pbinterno}.

The interior of $W$ may be empty or not empty; we consider the two cases separatedly. Let us start with $\mathring W =\emptyset$.
First, remark that
$\sigma$ cannot intersect~$W$, because otherwise (recall Proposition~\ref{local_omega}) 
$\sigma$ would intersect itself. As a consequence,
$\sigma$ is contained in only one connected component $P$ of~$S\setminus W$, and $W=\de P$.
Then, since $W$ does not reduce to a single pole and 
is connected (because $S$ is compact), we must have $W\cap S^0\neq\emptyset$.
Assume by contradiction that $W$ contains a closed geodesic~$\alpha$.
Since $\mathring W=\emptyset$, Proposition~\ref{local_omega}.(iii) implies that $\alpha$ is distinguished;
but then Lemma~\ref{distclosed} implies $W=\alpha$, impossible. So all the minimal sets contained in~$W$ are poles. 

By Proposition \ref{local_omega},
every point in $W\cap S^0$ must be contained in a unique maximal distinguished
geodesic contained in~$W$; moreover the $\omega$-limit sets of these
geodesic segments must be contained in~$W$, because $W$ is closed; we are going to prove that the $\omega$-limit sets of these geodesics
must be poles, completing the proof that $W$ is a boundary graph of saddle connections.

Let then $\rho$ be a distinguished geodesic and assume by contradiction that a point $w \in S^0$ belongs to the $\omega$-limit set $\omega (\rho)$ of $\rho$.
By the local form of the geodesics, and since $\omega(\rho)$, being contained in~$W$, has empty interior,
by Proposition \ref{local_omega} we find a unique geodesic segment
$\alpha$ through $w$ contained in $\omega (\rho)$. Take another geodesic segment transversal to $\alpha$ in~$w$.
Since $\rho$ must accumulate~$\alpha$, and $\rho$ is accumulated by~$\sigma$, it follows that
$\sigma$ must intersect the transversal arbitrarily near $w$ with both the directions of intersection.
This means that we have infinitely many disconnecting geodesic $2$-cycles,
with disjoint interiors and arbitrarily close to $\omega (\rho)$.
But by Remark~3.2
all these cycles must contain a pole, which is impossible.

Now let us consider the case $\mathring W \neq \emptyset$.
We know that $W$ is a closed union of leaves (and singular points) for $\Lambda$, and the same
is true for its boundary $\de W$, with $\de W \ne \emptyset$ because we are assuming $W \ne S$.
The minimal sets contained in~$\de W$ must be poles or closed geodesics, and $\de W$ cannot consist of poles only because otherwise we would have $W=S$. 
So, $\de W \cap S^0 \neq \emptyset$ and, by Proposition \ref{local_omega},
every $z \in \de W \cap S^0$ is contained in a unique geodesic segment 
contained in $\de W$. In particular, if by contradiction $\de W$ contained a closed geodesic 
we would have $W=\alpha$ (Lemma~\ref{distclosed}), impossible; so
arguing as above we see that each connected component of $\de W$ is a graph of saddle connections, that clearly 
cannot contain spikes (because $\sigma$ cannot cross~$\de W$), and we are done.
\ep
\medbreak

\proc{Remark 4.2.}
We have examples for all the cases of Theorem \ref{pbgeneral} (see Section~\ref{section_geod_toro}
 and \cite{ab}) except for cases \ref{pbtendcycle} and \ref{pbinterno}. Notice that in these cases
 if $p_0$ is a pole belonging to the boundary of the $\omega$-limit
 set of the geodesic~$\sigma$ then there must exists a neighbourhood $U$ of~$p_0$
 containing infinitely many distinct segments of~$\sigma$ (because $\sigma$ is
 accumulating~$p_0$ but it is not converging to it). The study of the local
 behavior of geodesics nearby a pole performed in~\cite{ab} implies that $p_0$
 must then be an irregular singularity or a Fuchsian singularity having
 $\Re\Res_{p_0}(\nabla)\ge-1$ (see \cite{ab} for the definitions, recalling
 that the connection we denote by $\nabla$ here is denoted by $\nabla^0$ in \cite{ab}).
Conversely, if all poles of~$\nabla$ are either apparent singularities or Fuchsian
singularities with real part of the residues less than $-1$ then cases \ref{pbtendcycle} and \ref{pbinterno} cannot happen.
 \medbreak
 
 \proc{Remark 4.3.}
There are examples of smooth line fields (and even vector fields) on Riemann
surfaces having integral curves accumulating graphs of saddle connections,
or having an $\omega$-limit set with not empty interior (but not covering everything); see, for example, \cite{st88}. 
\medbreak

\proc{Remark 4.4.}
When $S$ is the Riemann sphere, case \ref{pbinterno} cannot happen. Indeed,
assume by contradiction that a not self-intersecting geodesic $\sigma\subset S$
has an $\omega$-limit set $W$ with not-empty interior. Since $\sigma$ cannot
cross $\de W$, we can assume that $\sigma$ is contained in~$\mathring W$. Let
$w_0\in\de W$ be a smooth point of the boundary, and $\tau\subset W$ a
geodesic segment issuing from~$w_0$. The intersections between $\sigma$ and $\tau$ must
be dense in~$\tau$; then using suitable segments of $\sigma$ and $\tau$ we can construct
infinitely many disjoint geodesic 2-cycles. But in the Riemann sphere every 2-cycle is
disconnecting; thus by Remark~3.2 each of them must contain at least one pole. Using
the finiteness of the number of poles one readily gets a contradiction.
\medbreak

\section{$\omega$-limit sets and residues}\label{resd}

Corollary~\ref{ab4.2g} immediately gives a necessary condition for the existence of
geo\-desics having as $\omega$-limit set a disconnecting closed geodesic: 

\begin{cor}\label{sumresuno}
Let $S$ be a compact Riemann surface and $\nabla$ a meromorphic connection on $S$,
with poles $p_1, \dots, p_r \in S$, and set $S^0=S \setminus \{ p_1, \dots, p_r \}$.
Let $\sigma \colon [0, \eps_0) \to S^0$ be a maximal geodesic for $\nabla$, either
closed or having as $\omega$-limit set a closed geodesic. Assume that either $\sigma$ or
its $\omega$-limit set disconnects~$S$, and let $P$ be the part of $S$ containing~$\sigma$. Then
\begin{equation}\label{sumres}
\sum_{p_i \in P} \Re \Res_{p_i} (\nabla) = -1+ 2g_{\hat P}\;,
\end{equation}
where $\hat P$ is the filling of $P$.
\end{cor}

Furthermore, using Corollary~\ref{ab4.2g} as in \cite[Theorem~4.6]{ab} one
can easily get conditions on the residues that must be satisfied for the
existence of self-intersecting geodesics.
Our next aim is to find a similar necessary condition for the existence of
$\omega$-limit sets which are
boundary graphs of saddle connections. 

To do this,
we shall need a slightly more refined
notion of disconnecting graph, to avoid trivialities. To have an idea of the kind of situations we would like to avoid, consider
a boundary graph consisting in a non-disconnecting cycle and a small
(necessarily disconnecting) loop attached to a vertex. This boundary graph
is disconnecting only because of the trivial small loop, and if we shrink the
loop to the vertex we recover a non-disconnecting boundary graph.
To define a class of graphs where this cannot happen we first introduce a procedure that we shall call
\emph{desingularization of a graph $G$} (obtained as $\omega$-limit set of a geodesic) \emph{to a curve $\gamma$}.

Let $\sigma$ be a $\nabla$-geodesic without self-intersections and suppose that its $\omega$-limit set $G$
is a 
a boundary graph. If $G$ is a tree (that is, it does not contain cycles) then it does not disconnect $S$; so
let us assume that $G$ is not a tree (and thus it contains at least one cycle;
but this does not imply that $G$ disconnects~$S$). In particular, if $U$ is a small
enough open neighbourhood of~$G$ then $U\setminus G$ is not connected.

\proc{Remark 5.1.}
If a boundary graph contains a spike, it can be the boundary of only one connected open set.
\medbreak

The following Lemma describes the asymptotic behavior of $\sigma$ near $G$.

\begin{lemma}\label{yoga}
Let $S$ be a compact Riemann surface and $\nabla$ a meromorphic connection on $S$,
with poles $p_1, \dots, p_r \in S$. Let $S^0=S \setminus \{ p_1, \dots, p_r \}$.
Let $\sigma \colon [0, \eps_0) \to S^0$ be a maximal geodesic for $\nabla$.
Assume that $\sigma$ has no self-intersections and that its $\omega$-limit set $G$ is
a boundary graph containing at least one cycle. Let $U$
be a small connected open neighbourhood of $G$, with $U\setminus G$ disconnected and contained in~$S^0$.
Then the support of $\sigma$ is definitively contained in only one connected component $U_\sigma$ of $U \setminus G$.
\end{lemma}

\proc{Proof.}
If, by contradiction, $\sigma$ definitely intersects two different connected
components of~$U\setminus G$, then it must intersects $\de U$ infinitely many times;
being $\de U$ compact, $\sigma$ then would have a limit point in~$\de U$ which is disjoint from~$G$, impossible.
\ep
\medbreak

We shall denote by $U_\sigma$ the connected component of $U\setminus G$ given by
Lemma \ref{yoga}. In particular, notice that
we must have $G \subset \partial U_\sigma$, and that $U_\sigma$ must be contained in
the part~$P$ of~$S$ containing the support of~$\sigma$ whose boundary is~$G$.

Let us now describe the desingularization of $G$.
For every vertex $p_j$ of $G$, let
$B_j$ be a small open ball centered at $p_j$. Moreover, for each spike
$s_i\subseteq G$ let $C_i$ be a small connected open neighbourhood of $s_i$ with smooth boundary.
We can assume that all  the $B_j$ and $C_i$ are contained in an open
neighbourhood $U$ satisfying the hypotheses of Lemma~\ref{yoga}.
Then the union
of the following three sets is a Jordan curve~$\gamma$ in $S^0$, that we call a \emph{desingularization of $G$}:
\begin{itemize}
 \item $G \setminus \left( \bigcup_{j} B_j \cup \bigcup_i C_i\right)$;
 \item $\left(\bigcup_j \de B_j \cap U_\sigma\right) \setminus \bigcup_i C_i$;
 \item $\left(\bigcup_i \de C_i \cap U_\sigma\right) \setminus \bigcup_j B_j$.
\end{itemize}

The rationale behind this definition is the following: we take the graph and the boundary of the neighbourhoods
of the spikes outside the small balls at the vertices and we connect the pieces with small arcs
(which are contained in the boundaries of the small balls).
We see that, in particular, we can (uniformly) approximate the graph $G$ with desingularizing curves,
with respect to any global metric on the
compact Riemann surface $S$.

We are now ready to give the following definitions.

\proc{Definition 5.1.}
A boundary graph $G$ of saddle connections, which is the $\omega$-limit set of
a not self-intersecting $\nabla$-geodesic, is
\emph{essentially disconnecting} if every sufficiently close desingularization of $G$ disconnects $S$.
\medbreak

It is clear that, for sufficiently close desingularizations, if one of them is disconnecting then
all of them are. Because of this we may give the previous definition
without specifying which particular desingularization we are using. 

Furthermore, if $G$ is essentially disconnecting, then every sufficiently
close desingularization~$\gamma$ is the boundary of a well-defined part~$P_\gamma$ of~$S$,
the (closure of the) connected component of $S\setminus\gamma$ intersecting the part~$P$
whose boundary is~$G$ and containing the support of~$\sigma$. We can then consider the
filling $\hat P_\gamma$ obtained by glueing a disk along~$\gamma$. It is clear that, by
continuity, the genus of $\hat P_\gamma$ is the same for all close enough desingularizations of~$G$.

\proc{Definition 5.2.}\label{dgfP}
With a slight abuse of language we shall call \emph{genus of the filling of~$P$} this
common value, and we shall denote it by~$g_{\hat P}$. 
\medbreak

Before further investigating (essentially) disconnecting cycles, we give some examples and counterexamples of
possible boundary graphs, together with their desingularizations.

\proc{Example 5.1.}
Consider a single cycle $\gamma_1$ of saddle connections. Then, by Lemma~\ref{yoga},
 $\gamma_1$ locally disconnects $S$, and a geodesic $\sigma$ having $\gamma_1$ as
 $\omega$-limit set accumulates $\gamma_1$ locally on one side (say, staying inside
 an open set $U_1$). The definition of the desigularization of $\gamma_1$
 gives a curve $\omega_1$ in the same homology class of $\gamma_1$, contained in
 (the closure of) $U_1$. In particular,
 $\gamma_1$ disconnects $S$ if and only if $\omega_1$ does. So, for this cycle the
 properties of being disconnecting and of being essentially disconnecting are
 equivalent (see Lemma \ref{lemma_tutti_cicli_disconnettono} below for a more general result).
 \medbreak
 
 \proc{Example 5.2.}
 Let now $\gamma_2$ be the union of a non-disconnecting cycle $\gamma_2^a$ and a small,
 homotopically trivial, disconnecting cycle
 $\gamma_2^b$, intersecting only at a common pole. We see that $\gamma_2$ disconnects
 $S$ in two parts (the two given by $\gamma_2^b$). By
 Lemma \ref{yoga}, a geodesic $\sigma_2$ having $\gamma_2$ as $\omega$-limit
 set must be contained in only one of these components. Because of the fact that it must
 accumulate both
 $\gamma_2^a$ and $\gamma_2^b$, it cannot be contained in the small region bounded
 only by $\gamma_2^b$, and so it must be in the other one. So, we see that a
 desingularization of
 $\gamma_2$ is a curve $\omega_2$ contained in the closure of this last region,
 and in particular that $\gamma_2$ is not essentially disconnecting. 
\medbreak

\proc{Example 5.3.}
Let now $\gamma_3$ be the union of two non-disconnecting cycles $\gamma_3^a$
and $\gamma_3^b$ in the same
 homology class with exactly one vertex in common. We see that $\gamma_3$
 is disconnecting, and in particular that the fundamental group of one of the two components (that we call $U_3^1$)
 has rank one, while the fundamental group of the other component has the same
 rank as that of $\pi_1 (S)$ (in particular at least 2, since on the sphere we cannot find
 non-disconnecting cycles). A small neighbourhood of $\gamma_3$ is disconnected
 by it in three connected components, one contained in the component
 $U_3^1$ and two in the other component.
 The only one adherent to both $\gamma_3^a$ and $\gamma_3^b$ is the one contained in
 $U_3^1$, and so a geodesic $\sigma_3$ having $\gamma_3$
as $\omega$-limit set must live in $U_3^1$. So, we see that the construction of the desingularization
 of $\gamma_3$ gives a curve~$\omega_3$ contained in the closure of $U_3^1$.
 In particular, $\omega_3$ disconnects $S$, and so $\gamma_3$ is essentially
 disconnecting (which is coherent with the idea that small perturbations of
 $\gamma_3$ still disconnect $S$).
 \medbreak
 
 \proc{Example 5.4}
Finally, consider the cycle $\gamma_3$ of the previous example, and attach to it
at the intersection point, as in Example~4.2, a small homotopically trivial
disconnecting cycle~$\gamma_4^b$ not contained in~$U_3^1$. Lemma \ref{yoga} then implies that the graph $\gamma_4$
so obtained cannot be the $\omega$-limit set of any geodesic $\sigma_4$ on $S$. In fact, by the argument of the previous
 example, $\sigma_4$ should be contained in $U_3^1$, and thus it could not accumulate $\gamma_4^b$, contradiction.
\medbreak

\begin{figure}
 \centering
        \begin{subfigure}[b]{0.3\textwidth}
                \centering
                \includegraphics[width=\textwidth]{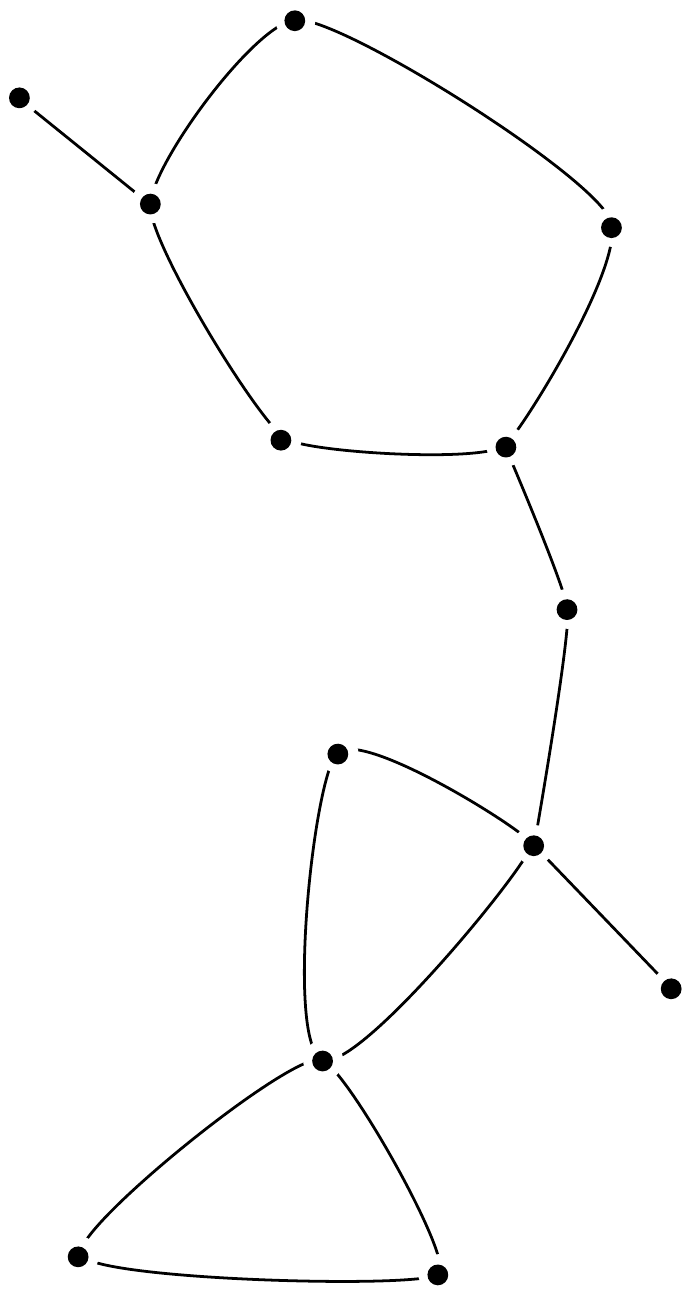}
                \caption{the graph}
                \label{fig:graph_empty}
        \end{subfigure}
        ~ 
        \begin{subfigure}[b]{0.3\textwidth}
                \centering
                \includegraphics[width=\textwidth]{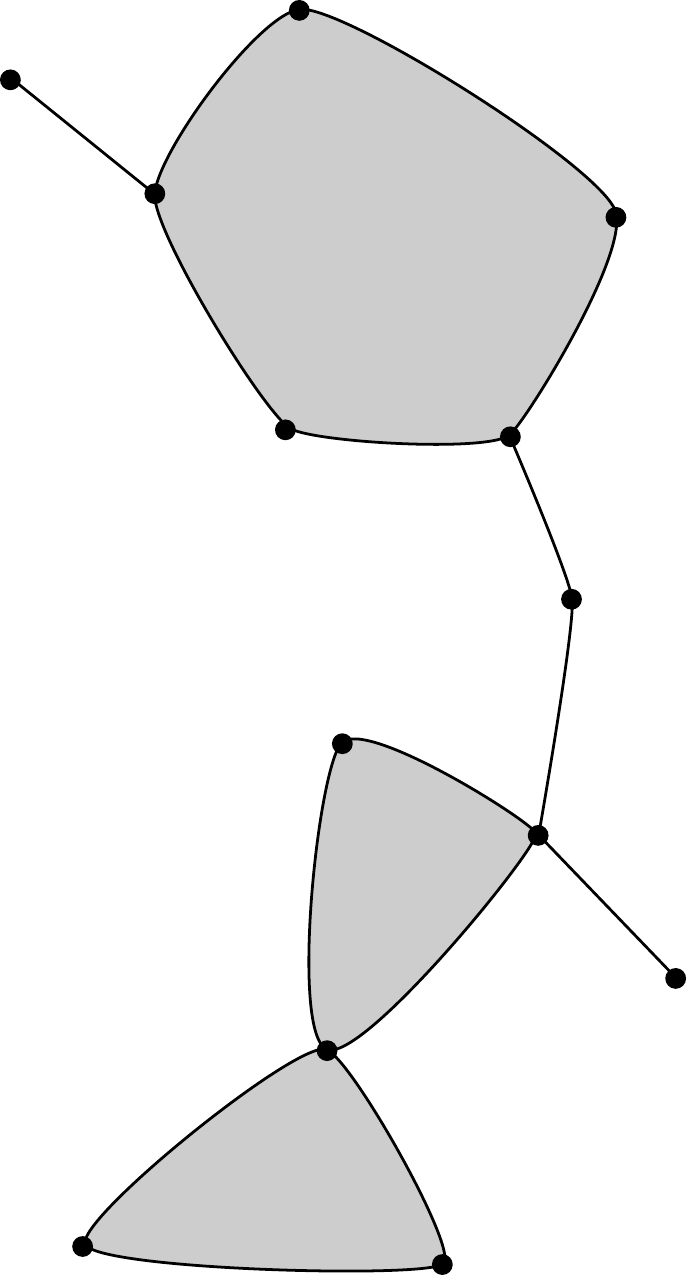}
                \caption{the black region}
                \label{fig:graph_full}
        \end{subfigure}
        ~ 
        \begin{subfigure}[b]{0.3\textwidth}
                \centering
                \includegraphics[width=\textwidth]{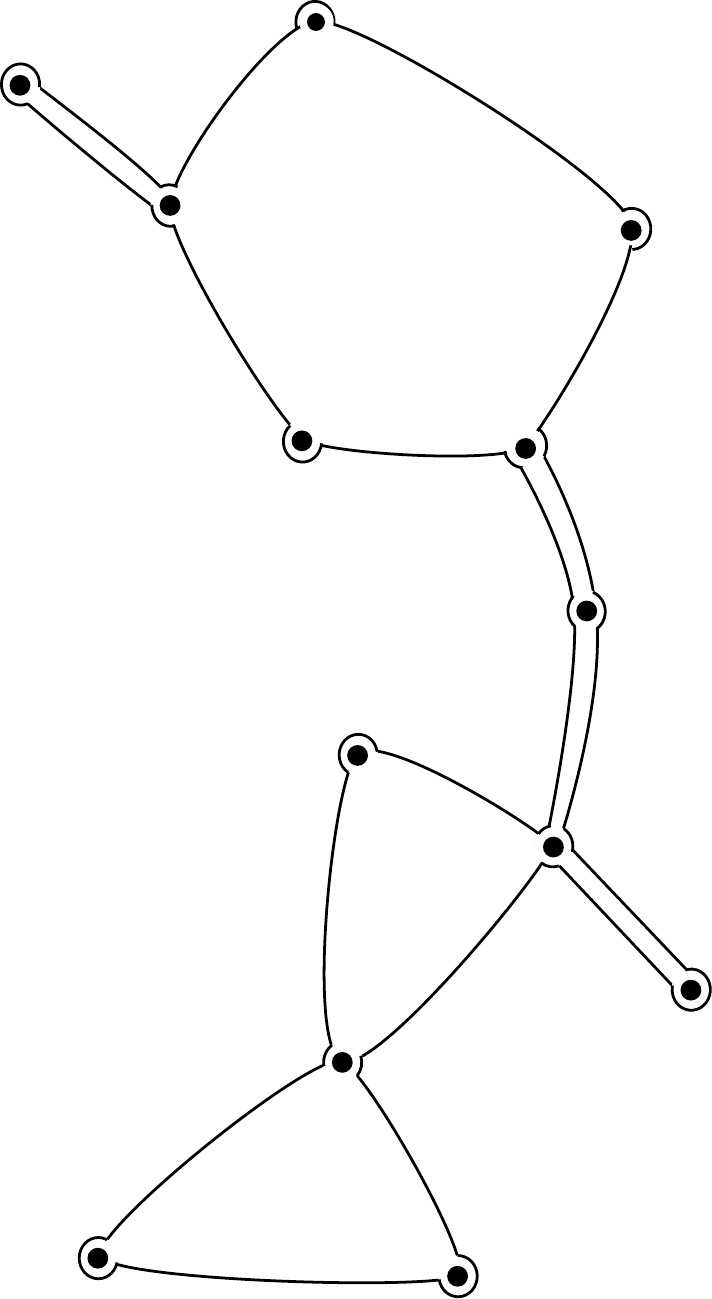}
                \caption{the curve $\omega$}
                \label{fig:graph_curve}
        \end{subfigure}
        \caption{Desingularization of a graph}\label{fig:desing_graph}
\end{figure}

The next lemma describes an important class of essentially disconnecting boundary graphs:

\begin{lemma}\label{lemma_tutti_cicli_disconnettono}
Let $G$ be a boundary graph of saddle connections in $S$, obtained
as $\omega$-limit set of a not self-intersecting geodesic, and such that every
cycle of $G$ disconnects~$S$. Then $G$ is essentially disconnecting.
\end{lemma}

\proc{Proof.}
Let $U_\sigma$ be as in Lemma~\ref{yoga}. 
Since, by assumption, every cycle in $G$ disconnects $S$, the complement
$S\setminus\Cc$ of any cycle $\Cc$ in $G$ has exactly two connected components,
only one disjoint from $U_\sigma$; let us color this component in black.

Furthermore, the fact that $G$ is a boundary graph implies that every edge
of $G$ is contained in at most one cycle;
as a consequence, $S\setminus G$ is formed by the part $P$
of $S$ containing $U_\sigma$ with $G$ as a boundary, and by the black connected components
determined by the cycles in $G$.

Now we
desingularize $G$ to a curve $\gamma$ in the following way, clearly
equivalent to the construction above:
near a pole connecting two (or more) cycles, we paint in black a little ball; and we replace
each spike by a little black strip following its path. In this way we get
a connected black region (see Figure~\ref{fig:desing_graph}).

Let $\gamma$ be the boundary of the black region we have constructed in this way.
Clearly $\gamma$ is (homotopic to) a desingularization of $G$ and,
because it separates a black region from a white one, by definition $\gamma$ disconnects $S$.
It follows that $G$ is essentially disconnecting.
\ep
\medbreak

So, we know that if every cycle is disconnecting then the graph is
essentially disconnecting, while (generalizing Example~4.2 above) we
can construct examples with an arbitrarily high number of disconnecting
cycles (and at least one non disconnecting cycle) which are not
essentially disconnecting.

We can now state the analogous of Equation \eqref{sumres} for case \ref{pbtendcycle}
of Theorem~\ref{pbgeneral}.

\begin{theo}\label{pbgeneral_res_grafi}
Let $S$ be a compact Riemann surface and $\nabla$ a meromorphic connection on $S$,
with poles $p_1, \dots, p_r \in S$. Let $S^0=S \setminus \{ p_1, \dots, p_r \}$.
Let $\sigma \colon [0, \eps_0) \to S^0$ be a maximal geodesic for $\nabla$ having
as $\omega$-limit set an essentially disconnecting boundary graph of saddle connections~$G$.
Let $P$ the part of $S$ having $G$ as boundary and containing the support of $\sigma$. Then
\begin{equation}\label{sumresbis}
\sum_{p_i \in P} \Re \Res_{p_i} (\nabla) = -1+ 2g_{\hat P}\;,
\end{equation}
where $g_{\hat P}$ is the genus of the filling of~$P$ (see Definition~5.2). 
\end{theo}

\proc{Proof.}
Let $U$ be a connected open neighborhood of the support of~$\sigma$ as in Lemma~\ref{yoga};
without loss of generality we can assume that the support of~$\sigma$ is completely contained in~$U_\sigma$.

Consider a generic point $z_0 \in G \cap S^0$. Using the local isometry $J$ near $z_0$ we see that we can find a
geodesic $\tau\colon  [0,\eps)\to U_\sigma$ issuing from $z_0$
intersecting $\sigma$ transversally infinitely many times and making an
angle~$\alpha$ close to~$\pi/2$ with the geodesic segment in~$G$
passing through~$z_0$.

Let $p_0\in U_\sigma$ be a point of intersection between $\sigma$ and $\tau$,
and let $p_1\in U_\sigma$ be the next (following~$\sigma$) point of
intersection between~$\sigma$ and~$\tau$; necessarily $p_1\ne p_0$ because
$\sigma$ has no self-intersections. The union of the segments of~$\sigma$
and~$\tau$ between $p_0$ and $p_1$ is a geodesic 2-cycle~$\gamma$ contained in~$U_\sigma$.
The first observation is that the two external angles of~$\gamma$ must have opposite signs.
Indeed, if they had the same sign (roughly speaking, if $\sigma$ leaves and comes back to~$\tau$ on
the same side of~$\tau$) then $\gamma$ would
be the boundary of a part of~$S$ contained in~$U_\sigma$, and thus
containing no poles, 
against Remark~3.2. The second observation is that $p_1$ must be closer (along~$\tau$)
to~$G$ than~$p_0$. Indeed, $\gamma$ disconnects $S$, because $G$ is
essentially disconnecting;
if $p_1$ were farther away from $G$ than $p_0$ then the support of~$\sigma$ after~$p_1$ would be contained
in the connected component of $S\setminus\gamma$ disjoint from~$G$,
impossible (notice that after $p_1$ the geodesic $\sigma$ cannot intersect $\gamma$ anymore,
because $\sigma$ has no self-intersections and any further intersection with the segment
of~$\tau$ contained between~$p_0$ and~$p_1$ would be on the same side of~$\tau$ and thus can be ruled out as before). 

Repeating this argument starting from $p_1$ we get a sequence $\{p_n=\sigma(t_n)=\tau(s_n)\}$ of points
of intersections between $\sigma$ and $\tau$, with $t_n$ increasing, $s_n$ decreasing,
and $p_n\to z_0$. Furthermore, if we denote by $\gamma_n$ the geodesic 2-cycle bounded by
the segments of~$\sigma$ and~$\tau$ between~$p_n$ and~$p_{n+1}$, we also know that the
external angles of~$\gamma_n$ have opposite signs. 

Since $G$ is essentially disconnecting, for $n$ large enough $\gamma_n$ disconnects~$S$;
let $P_n$ be the connected component of $S\setminus\gamma_n$ not containing~$G$. Notice that
the poles contained in~$P_n$ are exactly the same contained in~$P$, because 
$\gamma_n\subset U_\sigma\setminus G$. Furthermore, since the $\gamma_n$ are eventually
homotopic to desingularizations of~$G$, for $n$ large enough the genus of the filling of $P_n$
becomes equal to $g_{\hat P}$, the genus of the filling of~$P$. 

Up to a subsequence, we can assume that the directions of both $\sigma'(t_n)$ and $\sigma'(t_{n+1})$
converge to a direction in~$z_0$, necessarily the direction of the geodesic passing
through~$z_0$ contained in~$G$. Since the external angles of~$\gamma_n$ have
opposite signs, it follows that their sum must converge to~0. But, by
Corollary~\ref{ab4.3g} and the arguments above, this sum must be eventually constant;
so it must be~0, and \eqref{sumresbis} follows recalling again Corollary~\ref{ab4.3g}.
\ep
\medbreak

In particular, \cite[Theorem 4.6]{ab} is now a consequence of Theorems \ref{pbgeneral} and~\ref{pbgeneral_res_grafi}
(and Remark 4.4).

The formula \eqref{sumresbis} gives a necessary condition for
the existence of geodesics on $S$ having as $\omega$-limit set an essentially
disconnecting boundary graph of saddle connections.
In particular, if no sum of
real parts of residues is an odd integer
number of the form $-1+2g$
then no such geodesic can exist on~$S$.

The next statement contains a similar necessary condition for the existence
of isolated vertices in a boundary graph of saddle connections:

\begin{prop}\label{remark_residuo_spike}
Let $S$ be a compact Riemann surface and $\nabla$ a meromorphic connection on $S$.
Let $p$ be a pole belonging to a boundary graph of saddle connections~$G$ which is an
$\omega$-limit set for a not self-intersecting geodesic~$\sigma$. Suppose that $p$
is the vertex of only one arc~$\rho$ of the graph. Then
\[
\Re \Res_p (\nabla) = -\frac{1}{2}\;.
\]
\end{prop}

\proc{Proof.}
Let $\tau$ be a geodesic crossing transversally~$\rho$ at a point $q$ near $p$.
The geodesic~$\sigma$ must intersect both sides (with respect to~$\rho$) of~$\tau$ infinitely many times. Indeed,
if it would intersect eventually only one side, arguing as in the previous proof we would see that 
the external angles at consecutive intersections would have opposite signs; and this would force
the existence of points in~$G$ close to~$p$ but not belonging to~$\rho$,
impossible. 

So we can construct a sequence $\{\sigma_n\}$ of segments of~$\sigma$ connecting a point
$p_n$ on one side of~$\tau$ to a point $q_n$ on the other side of~$\tau$, with both
points of intersection converging to~$q$. Let $\gamma_n$ be the geodesic 2-cycle
consisting in~$\sigma_n$ and the segment of~$\tau$ between~$p_n$ and~$q_n$;
eventually, $\gamma_n$ must be contained in a simply connected neighbourhood of~$p$,
and thus it is the boundary of a part $P_n$ of~$S$ containing~$p$ and no other poles.
Notice that the genus of the filling of~$P_n$ is zero.

An argument similar to the one used above shows that the external angles of~$\gamma_n$
must have the same sign. Up to a subsequence, we can
also assume that the direction of $\sigma$ both at~$p_n$ and at~$q_n$ converge to the
same direction (which is the direction of~$\rho$ at~$q$); thus the sum of the external
angles of~$\gamma_n$ must tend to~$\pi$. Applying Corollary~\ref{ab4.3g} to~$P_n$ we get (that this sum must be constant and) that
\[
\pi = 2\pi \bigl( 1+ \Re\Res_p (\nabla) \bigr)\;,
\]
which gives the assertion.
\ep
\medbreak

So, in particular, if all the poles have the real part
of the residue different from $-1/2$, the graph cannot have vertices belonging to one arc only --- 
and thus in particular it cannot have spikes.

\proc{Remark 5.2.}
With the same argument as in Proposition \ref{remark_residuo_spike} we see
that if $p$ is a vertex belonging only to two arcs that are spikes, then the real part of its residue
must be zero.
\medbreak

\section{Non-compact Riemann surfaces}\label{non-comp}

In this section we shall generalize some of the 
results of Section~\ref{sezpbgeneral} to
non-compact Riemann surfaces. We shall limit ourselves to consider
Riemann surfaces $S$ that can be realized as an open subset of a compact Riemann surface~$\hat S$. In this setting,
all the results depending only on the local structure of geodesics still hold; for instance, 
Proposition~\ref{local_omega}, Lemmas~\ref{nointersect} and~\ref{distclosed}, and
Theorem~\ref{pblemma} go through without any changes.

It is also easy to adapt Theorem~\ref{prop_minimali_limite}. Indeed, the definition
of singular line field on a compact Riemann surface $\hat S$ used in \cite{hounie1981minimal}
is that of regular line field defined on an open subset~$S^0$ of~$\hat S$, and the singular
set is by definition $\hat S\setminus S^0$. In particular, if we consider
the $\omega$-limit set~$W$ of a geodesic $\sigma\colon I\to S^0$ in~$\hat S$
(and not only in~$S$) then $W\cap\de S$ is automatically invariant with respect to any line field associated to~$\sigma$.
Therefore arguing as in the proof of Theorem~\ref{prop_minimali_limite} we get
the following result.

\begin{theo}\label{prop_minimali_limite_nc}
 Let $S$ be an open connected subset of a compact Riemann surface $\hat S$.
Let $\nabla$ be a meromorphic connection on $S$, and denote by $S^0\subseteq S$
 the complement in~$S$ of the poles of $\nabla$. Let $\sigma \colon [0, \eps_0) \to S^0$ be
 a maximal not self-intersecting geodesic for $\nabla$.
 Then the possible minimal sets for $\sigma$ are the following:
 \begin{enumerate}
  \item a pole of $\nabla$;
  \item a closed simple curve contained in~$S^0$, which is a closed geodesic for $\nabla$,
  and in this case the $\omega$-limit set of $\sigma$ coincides with this closed geodesic;
  \item a point in $\de S\subset\hat S$;
  \item all of $\hat S$, and in this case $\hat S$ is a torus.
 \end{enumerate}
\end{theo}

Finally, in order to generalize Theorem~\ref{pbgeneral} we need to consider saddle
connections escaping to infinity, i.e., leaving any compact subset of~$S$.

\proc{Definition 6.1}
A \emph{diverging saddle connection} for a meromorphic connection on a Riemann surface $S$
is a maximal geodesic $\sigma\colon (-\eps_{-}, \eps_{+}) \to S^0$ satisfying one of the allowing conditions:
\begin{itemize}
\item[(a)] $\sigma (t)$
leaves every compact subset of $S$ for both $t \to -\eps_-$ and $t\to \eps_+$; or 
\item[(b)] $\sigma (t)$
leaves every compact subset of $S$ for $t \to -\eps_-$ and
tends to a pole for $t\to \eps_+$; or 
\item[(c)]  $\sigma(t)$ tends to a pole for $t\to -\eps_-$ and leaves
every compact for $t\to \eps_+$.
\end{itemize}
A \emph{diverging graph of saddle connections} is a connected planar
graph in $S$ whose vertices are poles and whose arcs are disjoint saddle connections,
with at least one of them diverging. A \emph{spike} is a saddle connection
of a graph which does not belong to any cycle of the graph.

A \emph{diverging boundary graph of saddle connections} (or \emph{diverging boundary graph})
is a diverging graph of saddle connections
which is also the boundary of a connected open
subset of $S$. A diverging boundary graph is \emph{disconnecting} if its complement in $S$ is not connected.
\medbreak

Using these definitions, the following theorem generalizes Theorem~\ref{pbgeneral} describing the possible $\omega$-limit sets
of geodesics for a meromorphic connection on a non-compact Riemann surface of this kind.

\begin{theo}\label{nocptpbgeneral}
Let $S$ be an open connected subset of a compact Riemann surface~$\hat S$.
Let $\nabla$ be a meromorphic connection on $S$, and denote by $S^0\subseteq S$
 the complement in~$S$ of the poles of $\nabla$.
Let $\sigma \colon [0, \eps_0) \to S^0$ be a maximal geodesic for $\nabla$.
Then one of the following possibilities occurs:
\begin{enumerate}
\item\label{nocptpbdiverges} $\sigma(t)$ definitely leaves every compact subset
of $S$ (i.e., it tends to the boundary of~$S$ in~$\hat S$); or
\item\label{nocptpbtendpole} $\sigma(t)$ tends to a pole in~$S$; or
\item\label{nocptpbclosed} $\sigma$ is closed; or
\item\label{nocptpbtendclosed} the $\omega$-limit set of $\sigma$ is the support of a closed geodesic; or
\item\label{nocptpbtendcycle} the $\omega$-limit set of $\sigma$ is a (possibly diverging) boundary graph of saddle connections; or
\item\label{nocptpbtutto} the $\omega$-limit set of $\sigma$ is all of $\bar{S}=S\cup\de S$; or
\item\label{nocptpbinterno} the $\omega$-limit set of $\sigma$ has non-empty interior and non-empty boundary in~$S$,
and each component of its boundary is a (possibly diverging) graph
of saddle connections with no spikes and at least one pole; or
\item\label{nocptpbinfint} $\sigma$ intersects itself infinitely many times.
\end{enumerate}
Furthermore, in cases \ref{nocptpbclosed} or \ref{nocptpbtendclosed} 
the support of $\sigma$ is contained in only one of the components
of the complement of the $\omega$-limit set, which is a part $P$ of $S$ having the $\omega$-limit set as boundary.
\end{theo}

\proc{Proof.}
Suppose that $\sigma$ is not closed (case \ref{nocptpbclosed}) and does not
self-intersect infinitely many times (case \ref{nocptpbinfint}).
Up to changing the starting
point, we can suppose that $\sigma \colon [0, \eps_0) \to S^0$ does not self-intersect at all.

If $\sigma$ definitely leaves every compact set of~$S$, we are in case~\ref{nocptpbdiverges}
and we are done. Otherwise, the $\omega$-limit set~$W$ of~$\sigma$ must intersect~$S$.
We can then argue as in the proof of Theorem~\ref{pbgeneral} using Theorems~\ref{pblemma}
and~\ref{prop_minimali_limite_nc}, and we are done.
\ep
\medbreak

\proc{Remark 6.1}
When the $\omega$-limit set is a non diverging graph
it is possible to find conditions on the residues of the poles involved analogous to
those described in Section~\ref{resd}, assuming it is disconnecting and the boundary of a relatively compact part. 
\medbreak

\section{Geodesics on the torus}\label{section_geod_toro}

In this section we study in detail the geodesics for holomorphic connections on the torus. We remark
that, by Theorem \ref{sum=-chi}, we cannot have holomorphic connections on Riemann surfaces different from the torus.
Moreover, Proposition \ref{prop_minimali_limite} tells us that this is the only case in which all the surface $S$ can be a minimal
set (and a minimal set for an associated line field, and thus Theorem~\ref{pbfoliaz} implies that the connection must be holomorphic).

So, our goals are: to characterize holomorphic connections on a torus, and to study their geodesics.

We recall that a complex torus can be realized as a quotient of $\C$ under the action of a rank-2 lattice, generated over $\R$
by two elements $\lambda_1$,~$ \lambda_2 \in \C$.
Without loss of generality we can assume 
that one of the two generators is 1 and the other is a $\lambda \in \C$
with $\Im \lambda >0$. We shall denote by $T_\lambda$ the torus associated to $\lambda$.

The next elementary result characterizes holomorphic connections on a torus.

\begin{lemma}
Every holomorphic connection on a torus is the projection of a holomorphic connection on the cover $\C$ represented
by a constant global form~$a\, dz$, for a suitable $a \in \C$.
\end{lemma}

\proc{Proof.}
Let $\nabla$ be a holomorphic connection on a torus~$T_\lambda$. 
Since the tangent bundle of a torus is trivial, we can find a global
holomorphic $(1,0)$-form $\eta$
representing $\nabla$ by setting $\nabla e =\eta\otimes e$, where $e$ is a global nowhere vanishing section of the tangent bundle.

Let $\pi\colon\C\to T_\lambda$ be the universal covering map, and let
$\tilde\eta=\pi^*\eta$; it is a holomorphic $(1,0)$-form on~$\C$. Thus
we can write $\tilde\eta= f\,dz$, for
a suitable entire function~$f$; we have to prove that $f$ is constant. 

Let $\tau_t\colon\C\to\C$ be given by $\tau_t(z)=z+t$. Since $\pi\circ\tau_1=\pi\circ\tau_\lambda=\pi$,
it follows that $\tau^*_1\tilde\eta=\tau^*_\lambda\tilde\eta=\tilde\eta$. But this is equivalent to 
$f(z+1)=f(z+\lambda)=f(z)$ for all~$z\in\C$, and thus $f$ must be
a bounded entire function, that is a constant. 
\ep
\medbreak

The next step consists in studying the geodesics for a holomorphic connection~$\nabla$ on a torus.
To do this, we shall study geodesics for the associated connection $\tilde{\nabla}$ on $\C$
represented by a constant form $a\, dz$ and then project
them to the torus (all geodesics on the torus are obtained in this way: see \cite[Proposition 3.1]{ab}).

Let $\tilde{\sigma}$ be a geodesic for
$\tilde{\nabla}$ issuing from a point $z_0 \in \C$ 
with tangent vector $v_0$.
We
first consider the trivial case $a=0$. In this case, the local isometry $J$ is given by
$J(z)= c z$, with $c \in \C^*$. So, applying equation \eqref{eqJ3} in Proposition~\ref{propJ}
we obtain
\[ 
c \tilde{\sigma}(t) = c v_0 t +cz_0\;,
\]
which means that $\tilde{\sigma}(t) = v_0 t +z_0$ and the geodesics are the euclidean ones.
Thus in this trivial case the geodesics on $T_\lambda$ are exactly projections of straight lines. 

If $a\ne 0$, the local isometry $J$ is given by $J(z)= \frac{1}{a} \exp(az)$.
We apply
again equation \eqref{eqJ3} in Proposition \ref{propJ}
to get
\[
\frac{1}{a} \exp(a \tilde{\sigma}(t)) = \exp (a z_0) v_0 t + \frac{1}{a} \exp (a z_0)\;,
\]
which we can solve to obtain
\begin{equation}\label{eqgeodpertoro}
\tilde{\sigma}(t)= z_0+\frac{1}{a} \log (1 + a v_0 t)\;,
\end{equation}
where $\log$ is 
the branch of the logarithm with $\log 1=0$ defined along the half-line $t \mapsto 1+av_0 t$.

So in general the geodesics for $\tilde{\nabla}$ are not straight lines, unless $av_0\in\R$.
However, every (projection of a) straight line in a torus is the support of
a geodesic for a suitable (non-trivial) holomorphic connection:

\begin{prop}\label{estoro}
Let $T_\lambda$ be a complex torus, and $\pi\colon\C\to T_\lambda$ the usual covering map.
Let $\ell$ be a straight 
line in~$\C$ issuing from $z_0\in\C$ and tangent to $v_0\in\C^*$. Then:
\begin{itemize}
\item $\pi(\ell)$ is the support of a geodesic $\sigma$ for the holomorphic
connection $\nabla$ represented by the global form $\eta=a\,dz$ with $a\in\C^*$
if and only if $v_0\in\R \bar{a}$;
\item $\pi(\ell)$ is the support of a closed geodesic $\sigma$ for the holomorphic
connection~$\nabla$ represented by the global form $\eta=a\,dz$ if and only if
$v_0\in\R \bar{a}$ and $\bar{a}\in\R^*(\Z\oplus\Z\lambda)$.
\end{itemize}
\end{prop}

\proc{Proof.}
For any $a\in\C^*$, the unique geodesic $\sigma$ in~$T_\lambda$ issuing from $\pi(z_0)$ and tangent
to $d\pi_{z_0}(v_0)$ is $\sigma=\pi\circ\tilde\sigma$, where $\tilde\sigma$ is given by \eqref{eqgeodpertoro}.
Then the support of~$\sigma$ coincides with $\pi(\ell)$ if and only if
$\tilde\sigma$ is a parametrization of~$\ell$, and this happens if and only if $av_0\in\R$,
which is equivalent to $v_0\in\R\bar{a}$.

Assume now $v_0\in\R\bar{a}$. Then $\sigma$ is closed if and only if $\tilde\sigma(t_0)=z_0+m+n\lambda$
for suitable $t_0\in\R^*$ and $m$,~$n\in\Z$, that is if and only if 
\[
a(m+n\lambda)=\log(1+av_0t_0)\;.
\]
Recalling that $av_0\in\R$, this can happens if and only if $\bar{a}=r(m+n\lambda)$ for a 
suitable $r\in\R^*$, and we are done. 
\ep
\medbreak

Some of the geodesics given by Proposition \ref{estoro} are closed, but none of them 
can be periodic.
We prove this in the next proposition, where we characterize closed geodesics.

\begin{prop}\label{torochiuse}
 Let $\nabla$ be a holomorphic connection on the torus
 $T_\lambda$  
 represented by a global form
 $\eta= a\, dz$.
 Let $\sigma \colon [0,\infty) \to T_\lambda$ be a (non-constant) closed geodesic
 for $\nabla$. Then:
 \begin{itemize}
  \item if $a\neq 0$ then the support of $\sigma$ is the projection of a
  straight line parallel to~$\bar{a}$ in the covering space~$\C$,
  and $\sigma$ cannot be periodic;
  \item if $a=0$ then $\sigma$ is periodic.
 \end{itemize}
\end{prop}

\proc{Proof.}
If $a=0$, the geodesics are the euclidean ones and so, once $\sigma$ is closed, it must also be periodic.

So let us study the problem when $a\neq 0$.
In order for $\sigma$ to be closed and non-trivial, \eqref{eqgeodpertoro} implies that we must have
\[
 \frac{1}{a} \log (1 + a v_0 \bar{t})= n+m\lambda
 \] 
for suitable $\bar{t} \in \R^*$
and $n$,~$m \in \Z$, not both zero, where $v_0=\tilde\sigma'(0)\ne 0$ and
$\tilde\sigma$ is the lifting of~$\sigma$ given by \eqref{eqgeodpertoro}.
Furthermore, since $\sigma$ is closed then $\tilde\sigma'(\bar{t})$ must be 
parallel to $\tilde{\sigma}'(0)$.
The derivative $\tilde{\sigma}'$ is given by
\begin{equation}\label{cp}
\tilde{\sigma}'(t)=  \frac{v_0}{1+a v_0 t}\;;
\end{equation}
so $1+a v_0 \bar{t}= e^{a(n+m\lambda)}$ must be real. Since $\bar{t}$ is real too, it follows that
$av_0$ must be real, and hence $\tilde\sigma$ is a parametrization of a straight line parallel to~$\bar a$, as claimed.

Finally, $\sigma$ is periodic if and only if $v_0=\tilde\sigma'(0)=\tilde\sigma'(\bar{t})$; but
by \eqref{cp} this can happen only if $a=0$, contradiction.
\ep
\medbreak

\proc{Remark 7.1.}
If the support of a geodesic $\sigma$ of a holomorphic connection on a torus is the projection of
a straight line, then $\sigma$ is either closed or everywhere dense.
In particular, Proposition~\ref{estoro} provides examples of cases \ref{pbclosed} and \ref{pbtutto} in Theorem~\ref{pbgeneral}. 
\medbreak

In the last theorem of this section we study the geodesics which are not
the projection of a straight line, and give
a complete description of the possible $\omega$-limit sets of geodesics for holomorphic connections on a torus.

\begin{theo}
 Let $\nabla$ be a holomorphic connection on the complex torus $T_\lambda$,
 represented by the global 1-form $\eta=a\, dz$. Let $\sigma\colon  [0,T) \to T_\lambda$
 be a maximal (non constant) geodesic for $\nabla$.
 \begin{itemize}
\item If $a \neq 0$ then the $\omega$-limit set of $\sigma$ is (the closure of) the projection on $T_\lambda$ of
a straight line $l\colon\R\to\C$ of the form $l(t)=\bar{a}t+b$ for a suitable $b\in\C$. In
particular:
\begin{enumerate}
 \item if $\bar{a} \in \R^*(\Z \oplus \lambda \Z) \setminus \{0\}$ then either
 \begin{itemize}
 \item $\sigma$ is closed, non periodic, and its support is the projection of the line~$l$; or
 \item $\sigma$ is not closed, but its $\omega$-limit set is a closed geodesic whose support is the projection
 of the line $l$;
 \end{itemize}
\item if $\bar{a} \notin \R^*(\Z \oplus \lambda \Z) \setminus \{0\}$ then
the $\omega$-limit set of $\sigma$ is all of $T_\lambda$.
\end{enumerate}

\item If $a=0$ then $\sigma$ is the projection of a line of the form $l(t)=v_0 t + b$, where $v_0 = \sigma' (0)$. In particular:
\begin{itemize}
 \item if $v_0 \in \R^+(\Z \oplus \lambda \Z)$ then $\sigma$ is closed and periodic;
 \item if $v_0 \notin \R^+ (\Z\oplus \lambda \Z)$ then $\sigma$ is not closed and its $\omega$-limit set is all
 of $T_\lambda$.
\end{itemize}
\end{itemize}

\end{theo}

\proc{Proof.}
If $a=0$, we know that the geodesics are the projection of euclidean straight lines, and the statement follows.
So let us assume $a\neq 0$.

Let $\tilde\sigma\colon[0,T)\to\C$ be the lifting of~$\sigma$ to the covering space~$\C$.
By \eqref{eqgeodpertoro} the lifting~$\tilde{\sigma}$ is given by
\[
\tilde{\sigma}(t)=z_0+ \frac{1}{a} \log (1 + a v_0 t)= z_0+\frac{\bar{a} \log(1+av_0 t)}{|a|^2}\;,
\]
where $v_0=\tilde\sigma'(0)$. When $t\to T$ it is easy to see
that $\Im\log(1+av_0t)=\mathrm{Arg}(1+av_0t)$ converges to a finite limit,
while $\Re\log(1+av_0t)=\log|1+av_0t|$ tends to $+\infty$ if $av_0\notin\R^-$, to $-\infty$ otherwise.
It follows that $\tilde\sigma$ is asymptotic to a line 
of the form $l(t)=\bar{a}t+b$ for a suitable $b\in\C$, and thus the $\omega$-limit
of~$\sigma$ is the closure of the projection of this line in~$T_\lambda$. The rest of the assertion
follows from Propositions~\ref{estoro} and~\ref{torochiuse}.
\ep


\begin{thebibliography}{[{\'E}ca81b]}

\bibitem[Aba01]{abate2001residual}
M.~Abate.
\newblock The residual index and the dynamics of holomorphic maps tangent to
  the identity.
\newblock {\em Duke Math. J.}, 107(1):173--207, 2001.

\bibitem[Aba10]{abasurvey}
M.~Abate.
\newblock Local discrete holomorphic dynamical systems.
\newblock In {\em Holomorphic dynamical systems}. G. Gentili, J. Guenot, G.
  Patrizio eds., Lecture Notes in Math., Springer, Berlin, vol. 1998, pp. 1--55, 2010.

\bibitem[ABT04]{abtann}
M. Abate, F. Bracci and F.~Tovena.
\newblock Index theorems for holomorphic self-maps.
\newblock {\em The Annals of Math.}, 159(2):819--864, 2004.

\bibitem[ABT08]{abtind}
M. Abate, F. Bracci and F.~Tovena.
\newblock Index theorems for holomorphic maps and foliations.
\newblock {\em Indiana Univ. Math. J.}, 57: 2999--3048, 2008.

\bibitem[AT03]{ab_parabolic}
M.~Abate and F.~Tovena.
\newblock Parabolic curves in {$\C^3$}.
\newblock {\em Abstr. Appl. Anal.}, 2003(5):275--294, 2003.

\bibitem[AT11]{ab}
M.~Abate and F.~Tovena.
\newblock {P}oincar{\'e}-{B}endixson theorems for meromorphic connections and
  homogeneous vector fields.
\newblock {\em J. Differential Equations}, 251:2612--2684, 2011.

\bibitem[AR14]{arra}
M. Arizzi and J. Raissy.
\newblock On {\'E}calle-HakimÕs theorems in holomorphic dynamics.
\newblock In {\em Frontiers in complex dynamics.} Eds. A. Bonifant, M. Lyubich, S. Sutherland.
Princeton University Press, Princeton, 2014, pp. 387--449.

\bibitem[Cam78]{camacho78}
C.~Camacho.
\newblock On the local structure of conformal mappings and holomorphic vector
  fields in {$\C^2$}.
\newblock {\em Ast{\'e}risque}, 59(60):83--94, 1978.

\bibitem[Cie12]{ciesielski2012}
K.~Ciesielski.
\newblock The {P}oincar{\'e}-{B}endixson {T}heorem: from {P}oincar{\'e} to the
  {X}{X}{I}st century.
\newblock {\em Cent. Eur. J. Math.}, 10(6):2110--2128,
  2012.

\bibitem[{\'E}ca81a]{ecalle1}
J.~{\'E}calle.
\newblock {\em Les fonctions r{\'e}surgentes. Tome I: Les alg{\`e}bres e
  fonctions r{\`e}surgentes}, volume 81-05.
\newblock Publ. Math. Orsay, 1981.

\bibitem[{\'E}ca81b]{ecalle2}
J.~{\'E}calle.
\newblock {\em Les fonctions r{\'e}surgentes. Tome II: Les fonctions
  r{\'e}surgentes appliqu{\'e}es {\`a} l'it{\'e}ration}, volume 81-06.
\newblock Publ. Math. Orsay, 1981.

\bibitem[{\'E}ca85]{ecalle3}
J.~{\'E}calle.
\newblock {\em Les fonctions r{\'e}surgentes. Tome III: L'{\'e}quation du pont et
  la classification analytique des objects locaux}, volume 85-05.
\newblock Publ. Math. Orsay, 1985.

\bibitem[Hak97]{hakim1997transformations}
M.~Hakim.
\newblock Transformations tangent to the identity. {S}table pieces of
  manifolds.
\newblock {\em Universit{\'e} de Paris-Sud, Preprint}, 145, 1997.

\bibitem[Hak98]{hakim1998analytic}
M.~Hakim.
\newblock Analytic transformations of {$({\C}^p,0)$} tangent to the identity.
\newblock {\em Duke mathematical journal}, 92(2):403--428, 1998.

\bibitem[Hou81]{hounie1981minimal}
J.~Hounie.
\newblock Minimal sets of families of vector fields on compact surfaces.
\newblock {\em J. Differential Geom.}, 16(4):739--744, 1981.

\bibitem[IY08]{ilyako}
Y.~Ilyashenko and S.~Yakovenko.
\newblock {\em Lectures on analytic differential equations}, volume~86.
\newblock American Mathematical Society, Providence, RI, 2008.

\bibitem[Mil06]{milnor}
J.~Milnor.
\newblock {\em Dynamics in one complex variable, 3rd ed.}
\newblock Princeton University Press, Princeton, 2006.

\bibitem[Mol09]{molino}
L. Molino.
\newblock The dynamics of maps tangent to the identity and with nonvanishing
  index.
\newblock {\em Trans. Amer. Math. Soc.},
  361(3):1597--1623, 2009.

\bibitem[Ron10]{rong}
F. Rong.
\newblock Absolutely isolated singularities of holomorphic maps of $\C^n$ tangent to the identity. 
\newblock {\em Pacific J. Math.}, 246(2):421--433, 2010.

\bibitem[Sch63]{pbgen}
A.~J. Schwartz.
\newblock A generalization of a {P}oincar\'e-{B}endixson theorem to closed
  two-dimensional manifolds.
\newblock {\em Amer. J. Math.}, 85(3):453--458, 1963.

\bibitem[Shc82]{shcherbakov1982topological}
A.A. Shcherbakov.
\newblock Topological classification of germs of conformal mappings with
  identical linear part.
\newblock {\em Moscow Univ. Math. Bull.}, 37:60--65, 1982.

\bibitem[ST88]{st88}
R.A.~Smith and E.S.~Thomas.
\newblock Transitive flows on two-dimensional manifolds.
\newblock {\em J. London Math. Soc. (2)}, 37:569--576, 1988.

\bibitem[Viv11]{vivas2011degenerate}
L.~Vivas.
\newblock Degenerate characteristic directions for maps tangent to the
  identity.
\newblock {\em Indiana Univ. Math. J.}, 61:2019-2040, 2012.

\end{thebibliography}
\end{document}